\documentclass[11pt, reqno]{amsart}
\usepackage{amsmath}
\usepackage{amssymb}

\numberwithin{equation}{section}
\begin{document}
\setlength{\baselineskip}{1.4em}

{\theoremstyle{plain}
    \newtheorem{thm}{\bf Theorem}[section]
    \newtheorem{pro}[thm]{\bf Proposition}
    \newtheorem{claim}[thm]{\sc Claim}
    \newtheorem{lemma}[thm]{\bf Lemma}
    \newtheorem{cor}[thm]{\bf Corollary}
    \newtheorem{Question}[thm]{\bf Question}
}
{\theoremstyle{remark}
    \newtheorem{remark}[thm]{\bf Remark}
    \newtheorem{example}[thm]{\bf Example}
}
{\theoremstyle{definition}
    \newtheorem{defn}[thm]{\bf Definition}
}


\def\To{\longrightarrow}
\def\height{\operatorname{ht}}
\def\reg{\operatorname{reg}}
\def\Hom{\operatorname{Hom}}
\def\Ext{\operatorname{Ext}}
\def\Proj{\operatorname{Proj}}
\def\BiProj{\operatorname{BiProj}}
\def\grade{\operatorname{grade}}
\def\Spec{\operatorname{Spec}}
\def\A{\mathcal A}
\def\Ann{\operatorname{Ann}}
\def\Ass{\operatorname{Ass}}
\def\rt{\operatorname{rt}}
\def\ord{\operatorname{ord}}
\def\img{\operatorname{im}}
\newcommand{\gothic}[1]{\frak #1}
\newcommand{\sfrac}[2]{\frac{\displaystyle #1}{\displaystyle #2}}
\newcommand{\smap}{\rightarrow\!\!\!\!\!\rightarrow}

\def\mm{{\frak m}}
\def\nn{{\frak n}}
\def\pp{{\frak p}}
\def\r{{\frak r}}
\def\O{{\mathcal O}}
\def\I{{\bf I}}
\def\J{{\bf J}}
\def\M{{\mathcal M}}
\def\L{{\mathcal L}}
\def\S{{\mathcal S}}
\def\NN{{\mathbb N}}
\def\PP{{\mathbb P}}
\def\ZZ{{\mathbb Z}}
\def\U{{\mathcal U}}
\def\T{{\bf T}}
\def\X{{\bf X}}
\def\8{\infty}
\def\xn{\x^{\n}}
\def\Xn{\T^{\n}}
\def\xm{\x^{\m}}
\def\e{{\bf e}}
\def\xni{\T^{\m+\e_i}}
\def\xnj{\T^{\m+\e_j}}
\def\n{{\bf n}}
\def\m{{\bf m}}
\def\k{{\bf k}}
\def\x{{\bf x}}
\def\y{{\bf y}}
\def\F{{\mathcal F}}
\def\length{\lambda}
\def\B{{\bf B}}
\def\a{{\bf A}}
\def\MM{{\bf M}}

\title{Homology multipliers and the relation type of parameter ideals}
\author{Ian M. Aberbach}
\address{Mathematics Department \\
    University of Missouri\\
    Columbia, MO 65211 USA}
\email{aberbach@math.missouri.edu}
\urladdr{http://www.math.missouri.edu/people/iaberbach.html}

\author{Laura Ghezzi}
\address{Mathematics Department \\
    University of Missouri\\
    Columbia, MO 65211 USA}
\email{ghezzi@math.missouri.edu}
\urladdr{http://www.math.missouri.edu/$\sim$ghezzi/}

\author{Huy T\`ai H\`a}
\address{Mathematics Department \\
    University of Missouri\\
    Columbia, MO 65211 USA}
\email{tai@math.missouri.edu}
\curraddr{Department of Mathematics, Tulane University, New Orleans LA 70118, USA. {\it Email:} {\tt tai@math.tulane.edu}. {\it URL:} {\tt http://www.math.tulane.edu/$\sim$tai/}}

\thanks{Part of this work was done when the authors visited MSRI. The authors would like to thank the institute for the financial support.}
\thanks{The first author was partially supported by  the National Security Agency and the University of Missouri Research Board.}

\begin{abstract} The relation type question,
raised by C. Huneke, asks whether for a complete equidimensional
local ring $R$ there exists a uniform number $N$ such that the
relation type of every ideal $I\subset R$ generated by a system of
parameters is at most $N$. Wang gave a positive answer to this question when the
non-Cohen-Macaulay locus of $R$ (denoted by ${\rm NCM}(R$)) has dimension zero. In this paper, we first present an example, due to the first author, which gives a negative answer to the question when $\dim {\rm NCM}(R) \ge 2$. The major part of our work is to investigate the remaining situation, i.e., when $\dim {\rm NCM}(R) =1$. We introduce the notion of homology multipliers and show that the question has a positive answer when $R/\A(R)$ is a domain, where $\A(R)$ is the ideal generated by all homology multipliers in $R$. In a more general context, we also discuss many interesting properties of homology multipliers.
\end{abstract}

\maketitle


\section{Introduction}\label{Intro}

Throughout this paper by ``ring'' we mean a commutative Noetherian ring with
identity.  

The existence of ``uniform bounds'' in Noetherian rings is an
interesting and important question. By uniform bounds we mean
statements which give some numerical bounds not just for one
ideal, but for all (or an infinite set of) ideals simultaneously.

In Noetherian rings we have an obvious finiteness condition, i.e.,
that every ideal is finitely generated;  there are, however, deeper forms
of finiteness which can be expressed in terms of uniform behavior.

Several types of uniform behavior have been demonstrated recently. See
for instance (this is by no means
a complete list) \cite{Hu, O1} (uniform Artin--Rees), \cite{Ra} (uniform
annihilation of local cohomology), \cite{Lai:1995, Wang:1997, Wang1} (uniform bounds on relation type).

In this paper we wish to extend the results of Lai and Wang 
concerning uniform bounds on relation type of parameter ideals in the
papers cited above. We are able to extend the
class of rings for which such uniform behavior exists.  Moreover, we
show that, in general, such uniform behavior should not be expected.

Let $R$ be a Noetherian  ring, and
let $I=(x_1,\dots,x_n)$ be an ideal of $R$. The Rees algebra
$R[It]$ of $I$ is a quotient of a polynomial ring over $R$. More
precisely, there is a canonical surjection $\phi:
R[T_1,\dots,T_n]\longrightarrow R[It]$ given by $T_i
\mapsto x_it$. By giving
degree $0$ to elements of $R$ and degree $1$ to $t$ and $T_i$ (for $1\leq i\leq n$) we have that $\phi$ is a homogeneous map, and so the kernel $Q$ of $\phi$ is a homogeneous ideal of $R[T_1, \dots, T_n]$.
The {\sl relation type} of $I$ is defined to be
$$\rt(I)=\min\{k|\  Q_k=Q\},$$ 
where $Q_k$ denotes the subideal of $Q$ generated by forms of
degrees $\leq k$. The relation type is independent of the choice of the generating set of $I$.

Let $F \in R[T_1, \dots, T_n]$ be a homogeneous form of degree $\delta$. It can be seen that
$$\phi(F) = F(x_1t, \dots, x_nt) = t^{\delta} F(x_1, \dots, x_n).$$
Thus, $F \in Q$ if and only if $F(x_1, \dots, x_n) = 0$. Therefore, by saying a {\sl relation on $x_1, \dots, x_n$} we mean a homogeneous form in $Q$. 

An ideal of relation type 1 is said to be of
linear type. Huneke \cite[Theorem~3.1]{Hu1} and Valla \cite[Theorem~3.15]{V}
proved that if $I$ is generated by a $d$-sequence, then $I$ is of
linear type. In particular, an ideal generated by a regular
sequence is of linear type.  Buchsbaum rings are precisely the rings
for which every parameter ideal is of linear type.

If $I$ is generated by a system of parameters (s.o.p.) in the local ring $R$ we say that $I$ is a {\sl parameter ideal}.
The following question was raised by C. Huneke.

\begin{Question}\label{q1}
{\rm (}The relation type question {\rm)} Let $R$ be a
complete equidimensional Noetherian ring of dimension $d$. Does
there exist an uniform number $N$ such that for every system of
parameters $x_1,\dots,x_d$ of $R$, $\rt(x_1,\dots,x_d)\leq N?$
\end{Question}

If a uniform bound as in Question~\ref{q1} exists,
we will say that $R$ satisfies {\sl bounded relation type}, or equivalently, $R$ has {\sl a uniform bound on relation type of parameter ideals}.

Question~\ref{q1} is closely connected to the strong uniform Artin-Rees
property.
Let $M\subseteq N$ be two finitely generated $R$-modules. The pair
$(M,N)$ is said to have the {\sl strong uniform Artin-Rees property} if
there exists an integer $k$ (depending on $M$ and $N$), such that
for all $R$-ideals $I$ and all $n\geq k$,
$$I^nM\cap N=I^{n-k}(I^kM\cap N).$$
There are several cases in which the strong uniform Artin-Rees
property holds (\cite{O1, DO, O2, Hu, Pl}), but Wang has shown in \cite{Wang1} that it does not
hold in general.  See \cite{Pl-V} for a recent summary and explication
of results relating to uniform Artin-Rees theorems.

Being a weaker version of the strong uniform Artin-Rees property
\cite{Lai:1995},  Question~\ref{q1} has attracted a great deal
of attention.

If $R$ is Cohen-Macaulay (CM) any system of parameters forms a regular sequence and so the relation type of any parameter ideal is 1. CM rings are characterized by the property that the local cohomology modules
$H_{\mathfrak m}^i(R)$ vanish for $i< \dim R$. The next step was 
to consider local rings $R$ such that $H_{\mathfrak m}^i(R)$ is
finitely generated (therefore of finite length) 
for all $i< \dim R$. 
Such rings are called
{\sl generalized Cohen-Macaulay} or {\sl rings with finite local cohomology} (f.l.c.).
Lai showed in \cite{Lai:1995} that bounded relation type 
holds for rings with finite local cohomology under the assumption
that the residue field is finite. In \cite{Wang1} Wang showed that
every $2$-dimensional Noetherian local
ring satisfies bounded relation type. 
Later he showed in \cite{Wang:1997} that bounded relation type holds
for rings with finite local cohomology without any restriction on
the residue field.

The first main result in this paper is to show that
bounded relation type does not hold in general. A
counterexample, due to the first author, has been known for some
time and is presented in Section~\ref{counterexample}. In this
example the non-CM locus of $R$ has dimension two
(and the counterexample easily generalizes to give counterexamples in
rings of arbitrarily high dimensional non-CM locus). On
the other hand, as observed at the beginning of Section \ref{bound-rt}, if $R$ is a complete equidimensional Noetherian ring, then $R$ has finite local cohomology if and only if $R$ has zero
dimensional non-CM locus.

The rest of the paper is devoted to studying the remaining case, i.e., when the 
non-CM locus of $R$ has dimension one. The methods of
\cite{Wang:1997} cannot be extended to this case. 
Wang uses strongly in his proof that if $R$ has f.l.c.~and $x\in R$ is a 
parameter, then $R/xR$ has f.l.c., and, moreover, the length of the
lower local cohomology modules in $R/xR$ can be bounded in terms of the lengths
of the lower local cohomology modules of $R$.  When $R$ is not generalized
CM then there is no uniform bound on the length of the local cohomology modules of
$R/xR$ as $x$ varies among parameters in $R$.
Hence, the starting
point of our work is an alternative proof of bounded relation type
 for rings with finite local cohomology which can be
generalized. We present this proof in Theorem~\ref{dparameters}.
We make use of ``homology multipliers'', defined in
Section~\ref{HM}, and of a ``Ramsey number'' combinatorial lemma, stated in
Section~\ref{ramsey-lemma}.

By a {\sl homology multiplier} we mean an element in
$R$ which annihilates all homology of complexes
satisfying the standard rank and height conditions. We denote by
$\A(R)$ the ideal of $R$ generated by homology multipliers. It
follows from a result of Hochster and Huneke that under mild conditions
on the ring $R$, $\A(R)$ is, up to radical,
 the defining ideal of the non-CM locus in
$R$ (see Corollary \ref{colonkillercor}). An important property of homology multipliers that we repeatedly use throughout the paper is the fact that, upon multiplying by a homology multiplier colons of monomial ideals in parameters behave as if the parameters were variables (see Remark \ref{expectedvalue}). 
We obtain several interesting results relating relation type and homology
multipliers.  We show that if an element of a s.o.p.~is ``adjusted'' by a
homology multiplier (and results in a new s.o.p.) then the relation type
is unchanged (see Theorem \ref{rt-equality}).  This result generalizes a result of Lai.  Also, the above process does not affect superficiality
(see Lemma \ref{super}).

We show in 
Section~\ref{counterexamplesec} that there are rings where $\dim(R/\A(R)) \ge 2$
which fail to satify bounded relation type. Hence  our attention is
focussed on the case  that $\dim R/\A(R) = 1$.

The new class of rings satisfying bounded relation type is given by
the theorem stated below, which appears in Section~\ref{sectionCM1} as
Theorem~\ref{CM1}.

\begin{thm} \label{intro-main}
Let $(R,\mm,k)$ be a formally unmixed local ring of dimension $d$
such that $\A(\hat{R})$ is a prime ideal of dimension one in the
completion $\hat{R}$ of $R$. Then $R$ has a uniform bound on
relation type of parameter ideals.
\end{thm}

We shall now briefly discuss the main steps of the proof. 
Our hypotheses allow us to assume that $R$ is complete.
Let
$I=(x_1,\dots,x_d)$ be a parameter ideal of  $R$, and let $\alpha$
be a homology multiplier. Theorem~\ref{rt-equality} shows that if
$x_1,x_2,\dots,x_d+\alpha$ is a system of parameters, then
$\rt(x_1,\dots,x_d)=\rt(x_1,x_2,\dots,x_d+\alpha).$ This
allows us to ``modify $x_d$ in a convenient way''. 
In particular, we replace $x_d$ by $w^t x'_d$ where $w$ is a reduction
of $\mm R/\A(R)$ and $x'_d$ is in a uniformly bounded small power of
$\mm R/\A(R)$.
 We do not have 
that $H_{\mm}^i(R)$ has finite length for all $i< d$ (as
in the generalized CM case), but since the ring
$R/\A(R)$ is a complete one-dimensional domain
  we can obtain uniform bounds for the lengths of
lower local cohomology modules of $R/x'_d R$.  

Given a relation $F(T_1,\ldots, T_d)$ on $x_1,\ldots, x_{d-1}, x_d = 
w^t x'_d$ we look at larger and larger ``partial sums'' and use
homology multipliers and uniformly bounded length of local cohomology
modules of $R/x_d'R$ to find a relation $G(T_1,\ldots, T_d)$ of uniformly
bounded degree such that the initial monomial of $G$ divides the initial
monomial of $F$.  Inductively, we obtain a bound (which does not
depend on the given system of parameters) on the highest degree
of a minimal generator of a relation on the system of parameters.

The general argument given is rather subtle, and certainly complicated.
However, the basic ideas in the general argument are already present 
in the argument for two parameters, in which case, the algorithm
is transparent.  We urge the reader to start with this case by 
reading Theorem~\ref{2parameters}.

In Corollary \ref{Fpure} we apply Theorem \ref{intro-main} to $F$-pure
rings.


\section{Bounded relation type does not hold in general}\label{counterexamplesec}

In this section, we will present an example showing that bounded relation type does not hold in general. In our example, the ring has non-CM locus of dimension 2, but it is possible to generalize to rings with non-CM locus of any dimension $\ge 2$.

\begin{example} \label{counterexample}
Let $R = k[[x,y,z,w]]$ where $w^2 = wz = 0$. Then, $R$ does not
have bounded relation type.
\end{example}

\begin{proof} Let $I_n = (x^{n-1}y+z^n, x^n, y^n)=(u_1,u_2,u_3)$. Clearly, $I_n$ is a parameter ideal for all $n \ge 1$.
 We will show that for every $n\in \mathbb{N}$, the ideal $I_n$ has relation type at least $n$.

We order monomials of $R[T_1, T_2, T_3]$ using lex and $T_1>T_2>T_3$.
 Consider the relation $wT_1^n-wT_2^{n-1}T_3$ on $u_1, u_2, u_3$
in $R$. Suppose that $wT_1^n-wT_2^{n-1}T_3$ can be written as a
combination of relations on $u_1, u_2, u_3$ of degree less than
$n$. Then there exists a relation $F(T_1,T_2,T_3)$ on $u_1, u_2,
u_3$ of the form
$$F(T_1,T_2,T_3)=wT_1^{n-1}+{\rm smaller \ terms}.$$
Let $\bar F$ be the image of $F$ in $(R/wR)[T_1,T_2,T_3]$. Since
$\bar F$ is a relation on $u_1, u_2, u_3$ in $R/wR$, and
$R/wR=k[[x,y,z]]$ is Cohen-Macaulay, we can write $\bar F=\bar
H_1\bar K_1+\bar H_2\bar K_2+\bar H_3\bar K_3$, where
$$K_1=-u_2T_1+u_1T_2,\  K_2=-u_3T_1+u_1T_3,\  K_3=-u_3T_2+u_2T_3$$ are
the Koszul relations on $u_1, u_2, u_3$ in $R$, $\bar K_1, \bar
K_2, \bar K_3$ are the Koszul relations on $u_1, u_2, u_3$ in
$R/wR$, and $\bar H_1, \bar H_2, \bar H_3$ are polynomials with
coefficients in $R/wR$. Hence $F-H_1K_1-H_2K_2-H_3K_3\in (w)R[T_1, T_2, T_3]$.
Write
$$F-H_1K_1-H_2K_2-H_3K_3=wF',$$ where $F'(T_1,T_2,T_3)\in R[T_1,T_2,T_3].$
Notice that $F'$ contains the term
$T_1^{n-1}$.

Since $wF'(u_1,u_2,u_3)=0$ in $R$, we have that
$F'(u_1,u_2,u_3)\in 0:_R w=(w,z)R$. Hence the image of $F'$ in
$(R/(w,z))[T_1,T_2,T_3]=k[[x,y]][T_1,T_2,T_3]$ is a relation on
$x^{n-1}y, x^n, y^n$ in $k[[x,y]]$. This implies that
$(x^{n-1}y)^{n-1}\in (x^n,y^n)(x^n,y^n,x^{n-1}y)^{n-2}$, a
contradiction (see also \cite[Example 6.1]{Wang1}).
\end{proof}

Note that the ring $R$ has an embedded prime of dimension two. It
would be of interest to find a ring with unbounded relation type
which is a domain.

\begin{remark} Notice that in Example \ref{counterexample}
the non-CM locus of $R$ is defined by $(z,w)$ and so
it has dimension 2.

In a similar fashion, we can construct examples of rings with
unbounded relation type and non-CM locus of any
dimension $\ge 2$. Let $S=k[[t_1,\dots,t_m,z,w]]$ where $w^2 = wz
= 0$ and $m\geq 2$. The proof of Example \ref{counterexample}
shows that $\rt (t_1^{n-1}t_2+z^n, t_1^n, t_2^n, t_3,\dots,t_m)\ge
n$ for every $n \ge 1$. The non-CM locus of $S$ is
defined by $(z,w)$ and so it has dimension $m$.
\end{remark}


\section{Homology multipliers and superficial sequences} \label{HM}

In this section, we introduce the notion of homology multipliers and investigate a number of properties that will be used later on. We shall also briefly recall the notion of filter-regular and superficial sequences. The proof for some of the results (e.g. Lemma \ref{resolution} and Corollary \ref{HMbasechange}) would be simpler if the ring contains a field. Our arguments work for rings of mixed characteristics as well.

Let $R$ be a Noetherian ring, and let $G_{\bullet}$ be a complex
of finitely generated free modules
$$G_{\bullet}:\ 0\rightarrow G_n\rightarrow G_{n-1}\rightarrow\dots\rightarrow
G_i\rightarrow \dots \rightarrow G_1\rightarrow G_0\rightarrow
0.$$ Denote by $\alpha_i$ the map from $G_i$ to $G_{i-1}$. Let
$b_i$ denote the rank of $G_i$, and let
$r_i=\sum_{t=i}^{n}(-1)^{t-i}b_i$ for $1\leq i\leq n$, while
$r_{n+1}=0$. Let $I_t(\alpha_i)$ be the ideal generated by the
$t\times t$ minors of $\alpha_i$. 

Recall that a complex $G_{\bullet}$ as above satisfies the {\sl
standard rank and height conditions} if rank $\alpha_i=r_i$ for
$1\leq i\leq n$, and height $I_{r_i}(\alpha_i) \geq i$ whenever
$1\leq i\leq n$ (see \cite{HH1:1990, HH2:1992}). For
simplicity, we denote $I_{r_i}(\alpha_i)$ by $I_i(G_\bullet)$. If $G_{\bullet}$ is acyclic we have that $I_i(G_\bullet)$ is an invariant of $H_0(G_{\bullet})$,
and localizes properly.

\begin{defn} Let $R$ be a Noetherian ring. We say that $z\in R$ is a {\sl homology multiplier} if for every
finite complex $G_{\bullet}$ satisfying the standard rank and
height conditions, $z$ annihilates the homology
$H_i(G_{\bullet})$, for every $i\geq 1$. We denote by $\A(R)$ the
ideal of $R$ generated by homology multipliers.
\end{defn}

The notion of homology multipliers is a generalization of Cohen-Macaulay multipliers (\cite{HH1:1990, HH2:1992}). Recall that
$R$ is {\sl equidimensional} if $\dim R/p=\dim R$ for all minimal
prime ideals of $R$. The following result is due to Hochster and Huneke. 

\begin{thm}\cite[Theorem 11.8]{HH1:1990} \label{colonkillerthm}
Let $R$ be an equidimensional local ring which is a homomorphic
image of a Gorenstein ring. Let $z$ be an element of $R$ such that
$R_z$ is Cohen-Macaulay. Then $z$ has a fixed power $z'$ such that
$z'$ is a homology multiplier.
\end{thm}

Homology multipliers are introduced to handle colons of monomial ideals in parameters (see Remark \ref{expectedvalue} for a precise statement). To see this, we first need the following lemma.

\begin{lemma}\label{resolution} Let
$R=\mathbb{Z}[Y_1,Y_1^{-1},\dots,Y_t,Y_t^{-1}][X_1,\dots,X_d]$,
where $Y_1,\dots Y_t,X_1\dots X_d$ are variables. Let $I$ be an
ideal generated by monomials in $X_1,\dots, X_d$.
\begin{enumerate}
\item There exists a resolution $F_\bullet$ of $R/I$ of length $\leq d$.
Furthermore, for every $1\leq i\leq d$,
\begin{align}
\prod_{1\leq j_1<\dots<j_i\leq d}(X_{j_1},\dots, X_{j_i}) \subseteq \sqrt{I_i(F_\bullet)}. \label{HMmember}
\end{align}
\item As a consequence, {\rm(\ref{HMmember})} holds for every
resolution $F_\bullet$ of $R/I$.
\end{enumerate}
\end{lemma}

\begin{proof} First we show the existence of a resolution of $R/I$
of length $\leq d$. This is clear if $d=1$, so we can assume that
$d>1$. We use induction on the number of generators $\mu(I)$ of
$I$. The claim is trivial if $\mu(I)=1$. Assume that $I=(I',m)$
where $m$ is a monomial with the largest possible exponent of
$X_d$, and $I'$ is a monomial ideal with $\mu(I')\geq 1$. Then
$I':m$ is a monomial ideal in $X_1,\dots X_{d-1}$. We have a
short exact sequence
$$0 \rightarrow R/m(I':m) \rightarrow R/I' \oplus R/mR \rightarrow
R/I\rightarrow 0. \label{eq4}$$ By induction on $\mu(I)$ we have
that $R/I' \oplus R/mR$ has a resolution of length $\leq d$. By
induction on $d$ we have that $R/m(I':m)$ has a resolution of
length $\leq d-1$. Hence the mapping cone is a resolution of $R/I$
of length $\leq d$. In particular this shows that if $G_{\bullet}$
is any resolution of $R/I$, then $I_{d+1}(G_\bullet)=R$.

Next we show that $(X_1,\dots,X_d)\subseteq
\sqrt{I_d(F_\bullet)}$. This is clear if $I_d(F_\bullet)=R$, so we
can assume that $I_d(F_\bullet)$ is a proper ideal. If there
exists $i$, with $1\leq i \leq d$, such that $X_i\notin
\sqrt{I_d(F_\bullet)}$, then localizing at $X_i$ we obtain a
contradiction to the existence of a resolution of length $\leq
d-1$ in the ring $\mathbb{Z}[Y_1,Y_1^{-1},\dots,Y_t,Y_t^{-1},
X_i,X_i^{-1}][X_1,\dots,X_{i-1},X_{i+1},\dots,X_d]$.

Now let $1\leq i \leq d-1$ and let $1\leq j_1<\dots<j_i\leq d-1$.
Localize at the $d-i$ variables $\{X_1,\dots,X_d\}$ $\setminus
\{X_{j_1},\dots, X_{j_i}\}$ and conclude as above. This finishes
the proof of ${\rm (1)}$.

${\rm (2)}$ follows form ${\rm (1)}$ since the ideals
$I_i(F_\bullet)$ do not depend on the resolution of $R/I$.
\end{proof}

\begin{cor} \label{HMbasechange}
Let $R$ be an equidimensional catenary local ring of dimension $d$
and let $x_1,\dots,x_d$ be a system of parameters of $R$. Let
$\phi : \ZZ[X_1, \dots, X_d] \longrightarrow R$ be the ring
homomorphism from a polynomial ring over $\ZZ$ to $R$ sending $X_i$ to $x_i$ for all $i$. Suppose $J$ is a monomial ideal in $\ZZ[X_1, \dots, X_d]$
and $G_\bullet$ is a resolution of $\ZZ[X_1, \dots, X_d]/J$ of length at most $d$. Let $F_\bullet = G_\bullet \otimes R$ be the complex obtained from $G_\bullet$ by base change. Then, $F_\bullet$ satisfies the standard rank and height
conditions.
\end{cor}

\begin{proof} By Lemma \ref{resolution}, for $1\leq i\leq d$,
$$\prod_{1\leq j_1<\dots<j_i \le
d}(X_{j_1},\dots, X_{j_i})\subseteq \sqrt{I_i(G_\bullet)}.$$ This
implies that $\prod_{1\leq j_1<\dots<j_i\leq d}(x_{j_1},\dots,
x_{j_i}) \subseteq \sqrt{I_i(F_\bullet)}.$ Since $F_{\bullet}$ is
obtained by tensoring $G_{\bullet}$ with $R$, the $r_i$'s are
unchanged. Also, since $R$ is equidimensional and catenary, we have
$${\rm height}\Big(\prod_{1\leq j_1<\dots<j_i\leq d}(x_{j_1},\dots, x_{j_i})\Big) \ge i.$$
Hence, $F_\bullet$ satisfies the standard rank and height conditions since $G_\bullet$ certainly does.
\end{proof}

The following remark will be used very often in the rest of the
paper.

\begin{remark} \label{expectedvalue}
Let $R$ be an equidimensional catenary local ring of dimension $d$
and let $x_1,\dots,x_d$ be a system of parameters of $R$. Let
$S=\ZZ[X_1, \dots, X_d]$ be a polynomial ring over $\ZZ$, and let
$\phi$ be as in Corollary \ref{HMbasechange}. For a d-uple of non negative integers $\n = (n_1, \dots, n_d)$, we denote by $\X^\n$ the monomial $X_1^{n_1} \dots X_d^{n_d}$, and by $\x^\n$ the monomial $x_1^{n_1} \dots x_d^{n_d}$. Let $I=(\X^{\n_1},\dots, \X^{\n_t})$ and let $\X^\m$ be any monomial in $S$ such that $\X^\m \not\in I$. By considering the minimal free resolution of $I + (\X^\m)$,
$$ \cdots \longrightarrow S^l \stackrel{\partial_2}{\longrightarrow} S^{t+1} \longrightarrow S \longrightarrow \sfrac{S}{I + (\X^\m)} \longrightarrow 0, $$
it can be seen that the colon ideal $I : (\X^\m)$ is generated by elements of the last row of the matrix of $\partial_2$. Thus, as a consequence of Corollary \ref{HMbasechange}, for any homology multiplier $z \in R$, we have
$$z(IR :_R \x^\m) = (I :_S \X^\m)R.$$
More generally, if $J = (\X^{\m_1}, \dots, \X^{\m_s})$ is another monomial ideal in $S$, and ${\mathfrak a}=IR=(\x^{\n_1},\dots,
\x^{\n_t})$ and ${\mathfrak b}=JR=(\x^{\m_1},\dots, \x^{\m_s})$, then $z({\mathfrak a}:_R {\mathfrak b})= (I:_S J)R$. Hence, up to multiplying by a homology multiplier, colons of monomials in $x_1, \dots, x_d$ behave as if the elements $x_1,\dots,x_d$ were variables.

\end{remark}

Let $R$ be an equidimensional catenary local ring with maximal
ideal $\mm$ and let $z\in \A(R)$. Let $x_1,\dots,x_k$ be a
sequence of elements of $\mm$ that is part of a system of
parameters of $R$. Then $z$ annihilates all the higher Koszul
homology $H_i(x_1,\dots,x_k,R), i\geq 1$. In particular,

$$z((x_1,\dots,x_{k-1})R:_Rx_kR)\subseteq (x_1,\dots,x_{k-1})R,$$
that is, $z$ is a Cohen-Macaulay multiplier (\cite{HH1:1990, HH2:1992}).

The next corollary shows that up to radical, $\A(R)$ is the
defining ideal of the non-Cohen-Macaulay locus in $R$.

\begin{cor}\label{colonkillercor}
Let $R$ be an equidimensional local ring which is a homomorphic
image of a Gorenstein ring. Let $z\in R$. Then $z\in
 \sqrt{\A(R)}$ if and only if $R_z$ is Cohen-Macaulay.
\end{cor}

\begin{proof} One direction is
Theorem \ref{colonkillerthm}. Conversely, let $z\in \sqrt{\A(R)}$,
and let $x_1/1,\dots,x_n/1$
 be a system of parameters of $R_p$, where $p$ is a prime ideal of
 $R$ such that $z\notin p$.
There exists a power $z'$ of $z$ such that
$z'((x_1,\dots,x_{k-1})R:_Rx_kR)\subseteq (x_1,\dots,x_{k-1})R$
for every $1\leq k\leq n$. Hence $x_1/1,\dots,x_n/1$ form a regular
sequence in $R_p$.
\end{proof}

We conclude this section recalling two definitions that we will
use later on.

\begin{defn} \label{filter_reg}
Let $S$ be a standard $\NN$-graded algebra over a local
ring $S_0$ (i.e., $S = S_0[S_1]$). The sequence of elements
$z_1,\ldots, z_n\in S$ is called {\sl filter-regular} if, for each
$i\ge 1$, $((z_1,\ldots, z_{i-1}):z_i)_n = (z_1,\ldots, z_{i-1})_n$ for $n \gg 0$.
\end{defn}

When $S_0$ has infinite residue field, then any homogeneous ideal
of $S$ may be generated by a filter-regular sequence (cf. \cite[Lemma 3.1]{Trung87}). If
$z_1,\ldots, z_n$ is a filter-regular sequence, then the Koszul
homology  modules $H_i(z_1,\ldots, z_n; S)$ vanish in sufficiently
high degree (cf. \cite [Lemma 4.7]{AbHu}).

We will be most interested in the case where $S$ is the associated
graded ring of an ideal generated by a system of parameters.

\begin{defn} \label{sup_high_deg}
Let $(R,\mm)$ be a local ring with infinite residue field and let
$I$ be an ideal. Recall that $x\in I$ is a {\sl superficial
element} for $I$ if for some integer $c$ and all $n \gg 0$,
$(I^n:x) \ \cap I^c = I^{n-1}$. Let $S = G(I)$, the associated graded ring of $I$. Let $x\in I$
and let $z = x +I^2 \in S_1$. Notice that $x$ is a superficial
element for $I$ if and only if $z$ is a filter-regular element.
Let $I=(x_1,\dots, x_d)$ and let $z_i = x_i +I^2 \in S_1$ for
$1\leq i\leq d$. If the sequence $z_1,\ldots, z_d$ is
filter-regular, we say that the sequence $x_1,\ldots, x_d$ is a
{\sl superficial sequence} for $I$.

\begin{remark} \label{sup-high-deg}
If $x_1,\ldots, x_d$ is a superficial sequence for $I$ then
there exists an integer $c'$ such that, if $1 \le i \le d$ and
$r_1 x_1 + \cdots + r_i x_i \in I^n$ with $n \ge c'$ and
$r_1,\ldots, r_i \in I^{c'}$ then we have $ r_1 x_1 + \cdots + r_i x_i =
r_1'x_1+\cdots + r_i' x_i$ with $r_1', \dots, r_i' \in I^{n-1}$. Notice that
the converse is not true in general.
\end{remark}

Given any set of generators $a_1,\ldots, a_d$ for $I$, there is a
Zariski--open set $U$ of $(R/\mm)^{d^2}$ such that setting $x_i =
\sum_{j=1}^d u_{ij} a_j$ with $(\overline{u_{ij}})_{1 \le i, j \le
d} \in U$, where $\overline{u_{ij}}$ is the image of $u_{ij}$ in $R/\mm$, gives a superficial sequence generating $I$.
\end{defn}


\section{An application of homology multipliers to relation type} \label{reltype_unchange}

In this section, we shall investigate two interesting properties of homology multipliers. First, we show that when an element of a system of parameters is changed by a homology multiplier (such that we still have a s.o.p.), the relation type is not changed (see Theorem \ref{rt-equality}). This provides a nice tool for studying relation type by modifying s.o.p.'s in a ``convenient'' way, which we will apply in Sections \ref{bound-rt} and \ref{sectionCM1}. Secondly, we prove that superficiality is also preserved by changing an element of a s.o.p. by a homology multiplier (see Lemma \ref{super}).

Throughout this section, $R$ is an equidimensional catenary local ring of dimension $d$ and $\langle x_1, \dots, x_d \rangle$ denotes a s.o.p. of $R$.

For a tuple $\n = (n_1, \dots, n_d)$, let $|\n | = n_1 + \dots +n_d$ and let $\x^\n$ denote
$x_1^{n_1} \cdots x_d^{n_d}$. From now on, we shall always use graded reverse lex monomial ordering. Suppose
$$ F(T_1, \dots, T_d) = \sum_{| \n | = n, \n \le \n_0} r_\n \T^\n $$
is a homogeneous form of degree $n$ in $R[T_1, \dots, T_d]$ with leading term $r_{\n_0} \T^{\n_0}$ which provides a relation on $x_1, \dots, x_d$, i.e. $F(x_1, \dots, x_d) = 0$ (if a monomial $\T^\n$ does not appear in $F$, we shall take $r_\n = 0$). Let $A$ be a new variable and write
$$ F(T_1, \dots, T_d) = F(T_1+A, T_2, \dots, T_d) - AG(T_1+A, T_2, \dots, T_d).$$

\begin{lemma} \label{map}
Let $R$ be an equidimensional catenary local ring of dimension $d$. Suppose $\x = \langle x_1, \dots, x_d \rangle$ and $\y = \langle y_1, x_2, \dots, x_d \rangle$ are s.o.p.'s such that $\alpha = y_1 - x_1 \in \A(R)$ is a homology multiplier. Then using the notation as above,
$$ \alpha G(y_1, x_2, \dots, x_d) \in ( \y^\m \ | \ \m < \n_0 ).$$
\end{lemma}

\begin{proof} Let $J_T = \big( \T^\n ~\big|~ |\n| = n, \n \le \n_0 \big) \subseteq \ZZ[T_1, \dots, T_d]$ and $J = \big( \x^\n ~\big|~ |\n| = n, \n \le \n_0 \big) = J_T \otimes R \subseteq R$. Then, $J_T$ is a stable ideal in the sense of \cite{EK} if we reverse the order of variables, i.e., if we list the variables as $T_d, \dots, T_1$, and so it admits an Eliahou-Kervaire graded free resolution. Notice that \cite{EK} in fact provided a graded free resolution of $J_T$ rather than of $\ZZ[T_1, \dots, T_d]/J_T$; there is, therefore, a shift of one index in our resolution compared to that given in \cite{EK}. By base change, $R/J$ admits the following complex
$$ \F_\bullet: \cdots \To F_2 \stackrel{\partial_2}{\To} F_1 \stackrel{\partial_1}{\To} R \To R/J \To 0, $$
where $F_i$'s are free $R$-modules, and $\partial_2$ is given by a matrix $M = M(F)$ where each column has exactly 2 non-zero entries and is of the form
\begin{eqnarray}
[0 \ \cdots \ -x_j \ \cdots \ x_i \ \cdots \ 0]^{\rm T} \ \text{for some} \ i < j. \label{col}
\end{eqnarray}

It follows from Corollary \ref{HMbasechange} that $\F_\bullet$ satisfies the standard rank and height conditions. Since $F(x_1, \dots, x_d) = 0$, we have $$[r_{\n_0} \ \cdots \ r_\n \ \cdots]^{\rm T} \in \ \text{ker} \ \partial_1,$$
where $[r_{\n_0} \ \cdots \ r_\n \ \cdots]^{\rm T}$ is the column vector of coefficients in $F$.
Since $\alpha$ annihilates $H_1(\F_\bullet)$, we have
$$ \alpha [r_{\n_0} \ \cdots \ r_\n \ \cdots]^{\rm T} \in \img \partial_2.$$
Let $K = K(F)$ be the number of columns of $M$ and let $C_1, \dots, C_K$ be the columns of $M$. We then can write
\begin{eqnarray}
\alpha [r_{\n_0} \ \cdots \ r_\n \ \cdots]^{\rm T} = b_1 C_1 + \dots + b_K C_K, \label{coefficient}
\end{eqnarray}
where $b_1, \dots, b_K \in R$.

Consider an arbitrary term $r_\n\T^\n = r_\n T_1^{n_1} \dots T_d^{n_d}$ of $F(T_1, \dots, T_d)$. We have
\begin{align*}
r_\n\T^\n & = r_\n(T_1+A-A)^{n_1} T_2^{n_2} \dots T_d^{n_d} \\
& = r_\n(T_1+A)^{n_1}T_2^{n_2} \dots T_d^{n_d} + r_\n T_2^{n_2} \dots T_d^{n_d}\big[\sum_{l=1}^{n_1} (-1)^l {n_1 \choose l}(T_1+A)^{n_1-l}A^l\big] \\
& = r_\n(T_1+A)^{n_1}T_2^{n_2} \dots T_d^{n_d} + Ar_\n \sfrac{\T^\n}{T_1^{n_1}} \big[\sum_{l=1}^{n_1} (-1)^l {n_1 \choose l}(T_1+A)^{n_1-l}A^{l-1}\big].
\end{align*}
Thus, we have
$$AG(T_1+A, T_2, \dots, T_d) = \sum_{r_\n \T^\n \ \text{is a term in} \ F(T_1, \dots, T_d)} Ar_\n \sfrac{\T^\n}{T_1^{n_1}} \big[\sum_{l=1}^{n_1} (-1)^l {n_1 \choose l}(T_1+A)^{n_1-l}A^{l-1}\big]. $$
This implies that
\begin{align}
\alpha G(y_1, x_2, \dots, x_d) = \sum_{r_\n \x^\n \ \text{is a term in} \ F(x_1, \dots, x_d)} \alpha r_\n \sfrac{\x^\n}{x_1^{n_1}} \big[\sum_{l=1}^{n_1} (-1)^l {n_1 \choose l}(x_1+\alpha)^{n_1-l}\alpha^{l-1}\big]. \label{contribution}
\end{align}
It follows from (\ref{coefficient}) that $\alpha r_\n \in (b_1, \dots, b_K)R$. Hence, $\alpha G(y_1, x_2, \dots, x_d)$ can be written as a combination of $b_1, \dots, b_K$ with coefficients in $R$.

Let $\e_l$ be the $l$-th unit vector of $\NN^d$ and suppose $C_1 = [0 \ \cdots \ -x_j \ \cdots \ x_i \ \cdots \ 0]^{\rm T}$ for some fixed $i < j$. Now, consider the contribution to $\alpha G(y_1, x_2, \dots, x_d)$ coming from $b_1C_1$ after substituting (\ref{coefficient}) to (\ref{contribution}). It follows from (\ref{coefficient}) that this contribution results from $r_{\m+\e_i} \x^{\m+\e_i} > r_{\m+\e_j} \x^{\m+\e_j}$ of $F(x_1, \dots, x_d)$, where $\e_i > \e_j$.
If $i > 1$, in which case $j > i > 1$, then in (\ref{contribution}), terms coming from $r_{\m+\e_i} \x^{\m+\e_i}$ and $r_{\m+\e_j} \x^{\m+\e_j}$ cancel each other since they have the same power of $x_1$. Suppose that $j > i = 1$. Substituting (\ref{coefficient}) to (\ref{contribution}), $r_{\m+\e_1} \x^{\m+\e_1} = r_{\m + \e_1} x_1^{m_1+1}x_2^{m_2} \dots x_d^{m_d}$ gives
$$-b_1 x_j\sfrac{\xm}{x_1^{m_1}} \big[ \sum_{l=1}^{m_1+1}(-1)^l{m_1+1 \choose l} (x_1+\alpha)^{m_1+1-l}\alpha^{l-1} \big], $$
and $r_{\m+\e_j}\x^{\m+\e_j}$ gives
$$b_1 x_1\sfrac{\xm x_j}{x_1^{m_1}} \big[ \sum_{l=1}^{m_1}(-1)^l{m_1 \choose l} (x_l+\alpha)^{m_1-l}\alpha^{l-1} \big]. $$
Thus, these 2 terms of $F(x_1, \dots, x_d)$ contribute to $\alpha G(y_1, x_2, \dots, x_d)$ the following:
$$b_1\sfrac{\xm x_j}{x_1^{m_1}} \big[ \underbrace{\sum_{l=1}^{m_1} (-1)^l {m_1 \choose l}(x_1+\alpha)^{m_1-l}\alpha^{l-1}x_1 - \sum_{l=1}^{m_1+1}(-1)^l {m_1+1 \choose l}(x_1+\alpha)^{m_1+1-l}\alpha^{l-1}}_{\displaystyle S} \big].$$
Note that ${m_1 + 1 \choose l} = {m_1 \choose l} + {m_1 \choose l-1}$. Therefore,
\begin{align*}
S & = \sum_{l=1}^{m_1} (-1)^l {m_1 \choose l}(x_1+\alpha)^{m_1-l}\alpha^{l-1}x_1 - \sum_{l=1}^{m_1+1}(-1)^l {m_1 \choose l}(x_1+\alpha)^{m_1+1-l}\alpha^{l-1} \\
& \quad - \sum_{l=1}^{m_1+1}(-1)^l {m_1 \choose l-1}(x_1+\alpha)^{m_1+1-l}\alpha^{l-1} \\
& = \sum_{l=1}^{m_1} (-1)^l {m_1 \choose l}(x_1+\alpha)^{m_1-l}\alpha^{l-1}x_1 - \sum_{l=1}^{m_1+1}(-1)^l {m_1 \choose l}(x_1+\alpha)^{m_1+1-l}\alpha^{l-1} \\
& \quad + \sum_{l=0}^{m_1}(-1)^l {m_1 \choose l}(x_1+\alpha)^{m_1-l}\alpha^l \\
& = \sum_{l=1}^{m_1}(-1)^l {m_1 \choose l} \big[ (x_1+\alpha)^{m_1-l}\alpha^{l-1}x_1 + (x_1+\alpha)^{m_1-l}\alpha^l \big] + (x_1+\alpha)^{m_1} \\
& \quad - \sum_{l=1}^{m_1+1}(-1)^l {m_1 \choose l}(x_1+\alpha)^{m_1+1-l}\alpha^{l-1}.
\end{align*}
We eventually have
\begin{align*}
S & = \sum_{l=1}^{m_1+1}(-1)^l {m_1 \choose l}(x_1+\alpha)^{m_1+1-l}\alpha^{l-1} - \sum_{l=1}^{m_1+1}(-1)^l {m_1 \choose l}(x_1+\alpha)^{m_1+1-l}\alpha^{l-1} + (x_1+\alpha)^{m_1} \\
& = (x_1+\alpha)^{m_1}.
\end{align*}
Hence, in $\alpha G(y_1, x_1, \dots, x_d)$, the term with factor $b_1$ is $b_1\y^{\m+\e_j}$. Since $\y^{\m+\e_j} < \y^{\n_0}$, we have $b_1 \y^{\m+\e_j} \in ( \y^\m \ | \ \m < \n_0 )$. A similar analysis works for the contribution to $\alpha G(y_1, x_2, \dots, x_d)$ coming from $b_2C_2, \dots, b_KC_K$. The lemma is proved.
\end{proof}

Suppose $\x = \langle x_1, \dots, x_d \rangle$ and $\y = \langle y_1, x_2, \dots, x_d \rangle$ are s.o.p.'s such that $y_1 - x_1 = \alpha \in \A(R)$. As before, let $F(T_1, \dots, T_d)$ be a homogeneous form of degree $N$ with leading term $r_{\n_0} \T^{\n_0}$ that gives a relation on $\x$. Let $M = M(F), K = K(F)$ and $C_1, \dots, C_K$ be as in Lemma \ref{map}. Suppose $C_{j_1}, \dots, C_{j_s}$ are the columns that contain $x_1$, i.e. for $i \le l \le s$, $C_{j_l}$ has the form
$$C_{j_l} = [0 \ \cdots \ -x_{jl} \ \cdots \ x_1 \ \cdots \ 0]^{\rm T}.$$
Suppose that in (\ref{coefficient}), $\alpha r_{\n_l}$, for each $l = 1, \dots, s$, lie in the same row as $x_1$. It follows from the proof of Lemma \ref{map} that, for each choice of $b_1, \dots, b_K$ in (\ref{coefficient}), we can write
$$ 0 = F(x_1, \dots, x_d) = F(y_1, x_2, \dots, x_d) + \sum_{l=1}^s b_{j_l} \y^{\n_l}.$$
Let $\tilde{F}_{b_1, \dots, b_K}(T_1, \dots, T_d) = F(T_1, \dots, T_d) + \sum_{l=1}^s b_{j_l} \T^{\n_l}$.
Then, for each choice of $b_1, \dots, b_K$ in (\ref{coefficient}), $\tilde{F}_{b_1, \dots, b_K}(T_1, \dots, T_d)$ provides a relation on $\y$. Notice that for any $b_1, \dots, b_K$, $\tilde{F}_{b_1, \dots, b_K}(T_1, \dots, T_d)$ and $F(T_1, \dots, T_d)$ always have the same leading term. We shall denote by $\Phi_{\alpha, b_1, \dots, b_K}$ the function which sends a relation $F(T_1, \dots, T_d)$ on $\x = \langle x_1, \dots, x_d \rangle$ to the relation $\tilde{F}_{b_1, \dots, b_K}(T_1, \dots, T_d)$ on $\y = \langle y_1, x_2, \dots, x_d \rangle$.

\begin{lemma}[cf. Theorem 3.6 of \cite{Lai:1995}] \label{property}
Let $\x = \langle x_1, \dots, x_d \rangle$ be a s.o.p.  and $\alpha \in \A(R)$. 
\begin{enumerate}
\item Suppose $F(T_1, \dots, T_d)$ is a relation on $\x$ and $1 \le l \le d$. Let $K = K(F)$ and $K' = K(T_lF)$. Then, for any choice of $b_1, \dots, b_K$ there exist $b_1', \dots, b_{K'}'$ such that
$$\Phi_{\alpha, b_1', \dots, b_{K'}'}(T_lF(T_1, \dots, T_d)) = T_l\Phi_{\alpha, b_1, \dots, b_K}(F(T_1, \dots, T_d)).$$
\item Suppose $F(T_1, \dots, T_d)$ and $H(T_1, \dots, T_d)$ are two relations on $\x$ of the same degree. Let $K' = K(F), K'' = K(H)$ and $K = K(F+H)$. Then, there exist choices of $b_1, \dots, b_K$, and $b_1', \dots, b_{K'}'$ and $b_1'', \dots, b_{K''}''$ such that
$$\Phi_{\alpha, b_1, \dots, b_K}(F+H) = \Phi_{\alpha, b_1', \dots, b_{K'}'}(F) + \Phi_{\alpha, b_1'', \dots, b_{K''}''}(H).$$
\item Suppose $F(T_1, \dots, T_d)$ is a relation on $\x$ and $P(T_1, \dots, T_d)$ is any polynomial. Let $K = K(F)$ and $K' = K(PF)$. Then, there exist choices of $b_1, \dots, b_{K'}$ and $c_1, \dots, c_{K}$ such that
$$\Phi_{\alpha, b_1, \dots, b_{K'}}(P(T_1, \dots, T_d)F(T_1, \dots, T_d)) = P(T_1, \dots, T_d)\Phi_{\alpha, c_1, \dots, c_{K}}(F(T_1, \dots, T_d)).$$
\end{enumerate}
\end{lemma}

\begin{proof} It is easy to see that (3) is a consequence of (1) and (2). We shall first prove (1). Suppose 
$$F(T_1, \dots, T_d) = \sum_{\n \le \n_0} r_\n \T^\n,$$
with leading monomial $\T^{\n_0}$. Let $M = M(F)$ and $M' = M(T_lF)$ be the presentation matrices associated to relations $F$ and $T_lF$, respectively, as obtained in Lemma \ref{map}. Suppose $C = [0 \ \dots \ -x_j \ \dots \ x_i \ \dots \ 0]^{\rm T}$ is a column of $M$ whose entries $-x_j$ and $x_i$ give the relation between monomials $\x^{\n_1}$ and $\x^{\n_2}$ of $F(x_1, \dots, x_d)$. Then, in $M'$, there is a corresponding column $C' = [0 \ \dots \ -x_j \ \dots \ x_i \ \dots \ 0]^{\rm T}$ (with more $0$'s) whose $-x_j$ and $x_i$ entries give the relation between monomials $\x^{\n_1 + \e_l} = x_l \x^{\n_1}$ and $\x^{\n_2+\e_l} = x_l\x^{\n_2}$ of $x_lF(x_1, \dots, x_d)$. By re-indexing, if necessary, we may assume that $C_1', \dots, C_K'$ are columns in $M'$ corresponding to columns $C_1, \dots, C_K$ of $M$. Now, in the presentation obtained from the relation $T_lF(T_1, \dots, T_d)$ similar to (\ref{coefficient}), we may pick $b_i' = b_i$ for $i = 1, \dots, K$ and $b_i' = 0$ for $i = K+1, \dots, K'$. (1) then follows from the construction of $\Phi$. 

It remains to prove (2). Without loss of generality, we may assume that $K' \ge K''$. If the leading terms of $F$ and $H$ do not cancel each other (which implies $K = K'$), then (2) follows from the construction of functions $\Phi_{\alpha, b_1, \dots, b_K}$ by taking the tuple
$$(b_1, \dots, b_K) = (b_1', \dots, b_{K'}') + (\underbrace{0, \dots, 0}_{K' - K''}, b_1'', \dots, b_{K''}'').$$
Suppose now that $F(T_1, \dots, T_d) = \sum_{\m \le \n_0}r_\m \T^\m$ and $H = \sum_{\m \le \n_0}r_\m'\T^\m$, and  $r_{\n_1}\T^{\n_1}$ and $r_{\n_1}'\T^{\n_1}$ are the highest terms in $F$ and $H$ that do not cancel (with $\n_1 < \n_0$). In this case, $K < K' = K''$ To prove (2), we only need to show that there are choices of $b_1, \dots, b_K$ and $b_1', \dots, b_{K'}'$ and $b_1'', \dots, b_{K''}''$ such that the leading term of $\Phi_{\alpha, b_1', \dots, b_{K'}'}(F) + \Phi_{\alpha, b_1'', \dots, b_{K''}''}(H)$ is the same as that of $\Phi_{\alpha, b_1, \dots, b_K}(F+H)$, which is $(r_{\n_1}+r_{\n_1}')\T^{\n_1}$.

Let $C_1, \dots, C_{K'}$ ($K' = K''$) be the columns of the matrix of $\partial_2$'s corresponding to $F$ and $H$ as in Lemma \ref{map}. Suppose $b_1', \dots, b_{K'}'$ and $b_1'', \dots, b_{K''}''$ are the coefficients in equalities of the form (\ref{coefficient}) corresponding to $F$ and $H$. We have
\begin{align}
\left( \begin{array}{c} 0 \\ \vdots \\ \alpha (r_{\n_0}+r_{\n_0}') \\ \vdots \\ \alpha (r_{\n_1}+r_{\n_1}') \\ \vdots \end{array} \right) = (b_1'+b_1'') C_1 + \dots + (b_{K'}'+b_{K''}'') C_{K'}. \label{HMprop}
\end{align}
Suppose in (\ref{HMprop}), $C_1, \dots, C_L$ are columns that have at least a non-zero entry higher than $\alpha(r_{\n_1}+r_{\n_1}')$. Then, $b_1' + b_1'' = \dots = b_L' + b_L'' = 0$. It can be easily seen that $K = K' - L$. We can now pick $b_1 = b_{L+1}' + b_{L+1}'', \dots, b_K = b_{K'}' + b_{K''}''$ and the required equality follows from (\ref{HMprop}). The lemma is proved.
\end{proof}

\begin{thm}[cf. Theorem 3.1 of \cite{Lai:1995}] \label{rt-equality}
Let $(R,\mm,k)$ be a local ring.  Assume that
$\x = \langle x_1,\dots, x_d \rangle$ and $\y = \langle y_1,x_2,\dots, x_d \rangle$ are both s.o.p.'s such that $\alpha = y_1 - x_1 \in \A(R)$.  Then
$$\rt(x_1,\dots, x_d) = \rt(y_1,x_2,\dots, x_d). $$
\end{thm}

\begin{proof} It is enough to show that of $\rt(x_1, \dots, x_d) \le \rt(y_1, x_2, \dots, x_d)$ (since we then can apply the inequality for $-\alpha$). Suppose $r = \rt(y_1, x_2, \dots, x_d)$. Let $F(T_1,
\dots, T_d)$ be a relation on $x_1, \dots, x_d$ of degree $N > r$. Let $K = K(F)$ and let $\tilde{F}_{b_1, \dots, b_K}(T_1, \dots, T_d) = \Phi_{\alpha, b_1, \ldots, b_K}(F(T_1, \dots, T_d))$ for some choice of $b_1, \dots, b_K$. Then, as shown before, $\tilde{F}_{b_1, \dots, b_K}$ gives a relation on $\y = (y_1, x_2, \dots, x_d)$. Since $\tilde{F}_{b_1, \dots, b_K}$ has degree $N > r$, it can be written as
$$ \tilde{F}_{b_1, \dots, b_K}(T_1, \dots, T_d) = \sum_{i=1}^m P_i(T_1, \dots, T_d) H_i(T_1, \dots, T_d), $$
where $P_i$'s are polynomial in $T_1, \dots, T_d$, and $H_i$'s provide relations on $\y = (y_1, x_2, \dots, x_d)$ with $\deg H_i \le r$ for all $i$. For $i = 1, \dots, m$, let $K_i = K(H_i)$. It follows from Lemma \ref{property} that there exist choices of $c_1, \dots, c_K$ and $b_{i1}, \dots, b_{iK_i}$ such that
$$ \Phi_{-\alpha, c_1, \dots, c_K}(\tilde{F}_{b_1, \dots, b_K}) = \sum_{i=1}^m P_i(T_1, \dots, T_d) \Phi_{-\alpha, b_{i1}, \dots, b_{iK_i}}(H_i(T_1, \dots, T_d)). $$
By the definition, $\Phi_{-\alpha, b_{i1}, \dots, b_{iK_i}}(H_i(T_1, \dots, T_d))$ gives a relation on $x_1, \dots, x_d$ for each $i$, and
$$\deg \Phi_{-\alpha, b_{i1}, \dots, b_{iK_i}}(H_i(T_1, \dots, T_d)) \le r.$$ Moreover, since $\tilde{F}_{b_1, \dots, b_K}$ and $F$ have the same leading term, $\Phi_{-\alpha, c_1, \dots, c_K}(\tilde{F}_{b_1, \dots, b_K})$ and $F$ also have the same leading term. Thus, we can write
$$F(T_1, \dots, T_d) = \sum_{i=1}^m P_i(T_1, \dots, T_d) \Phi_{-\alpha, b_{i1}, \dots, b_{iK_i}}(H_i(T_1, \dots, T_d)) + F'(T_1, \dots, T_d), $$
where $F'(T_1, \dots, T_d)$ gives a relation on $x_1, \dots, x_d$, and has a smaller leading term than that of $F(T_1, \dots, T_d)$. Repeating this process, we eventually will get to the situation when $F'(T_1, \dots, T_d) = 0$, i.e. $F(T_1, \dots, T_d)$ is a combination of relations on $x_1, \dots, x_d$ with degrees at most $r$. Hence, $\rt(x_1, \dots, x_d) \le r = \rt(y_1, x_2, \dots, x_d)$. The theorem is proved.
\end{proof}

We now refer the reader to Definition \ref{sup_high_deg} for the definition of superficial sequence. The following lemma will be useful later on.

\begin{lemma} \label{super}
Let $\x = \langle x_1, \dots, x_d \rangle$ be a s.o.p.  such that $x_d$ is
superficial for $(\x)$. Suppose $\y = \langle x_1, \dots, x_{d-1}, y_d \rangle$ is
a s.o.p.  such that $\alpha = y_d - x_d \in \A(R)$. Then $y_d$ is
superficial for $(\y)$.
\end{lemma}

\begin{proof} 
Suppose $c \in \NN$ is an integer such that $[(\x)^n:x_d] \cap (\x)^c = (\x)^{n-1}$ for all $n > c$.  Also, let $k\in \NN$ be an integer
given by the Artin-Rees lemma for the modules $x_dR \subseteq R$ and
the ideal $(x_1,\ldots, x_{d-1})$, i.e., for $m \ge k$ we have
$$(x_1, \dots, x_{d-1})^m \cap (x_d) \subseteq x_d(x_1, \dots, x_{d-1})^{m-1}.$$  
We first observe that $(0:x_d) \subseteq H^0_\mm(R)$ since if $x_du =0$
we have $x_d(\x)^cu = 0 \subseteq (\x)^n$ for all $n>c$.  By
superficiality, $(\x)^cu \in (\x)^{n-1}$ for all $n \gg 0$, so by
the Krull intersection theorem, $(\x)^cu =0$.  Since $(\x)$ is $\mm$-primary,
$u \in H^0_\mm(R)$.  We let $t \in \NN$ be such that $\mm^t \cap H^0_\mm(R)
= 0$.

We will now show that if $s \in (\y)^{c + k + t}$   and
$y_d s \in (\y)^m$ for $m > c + k + t$  then $s \in (\y)^{m-1}$.
Since we then have $y_d(s -s') \in (x_1, \dots, x_{d-1})^m$ where
$s' \in (\y)^{m-1}$, it suffices to assume that 
$s \in (x_1, \dots, x_{d-1})^m:y_d$.
 Since $\alpha$ is a homology multiplier, we then
have $\alpha s \in (x_1, \dots, x_{d-1})^m$. Thus,
\begin{align}
sx_d = sy_d - \alpha s \in (x_1, \dots, x_{d-1})^m. \label{eq361}
\end{align}

From the Artin-Rees lemma we see that $sx_d \in (x_1, \dots, x_{d-1})^m \cap (x_d) \subseteq x_d (x_1, \dots, x_{d-1})^{m-k}$, hence
$s \in [(x_1, \dots, x_{d-1})^{m-k} + (0:x_d)] \cap (\y)^{t}$.  If
we write $s = s_1 + s_2$ where $s_2 \in (0:x_d)$, then $s_2 \in
(0:x_d) \cap \mm^t \subseteq H^0_\mm(R) \cap \mm^t= 0$ 
(since $m-k \ge t$).  Thus $s \in (\x)^c$,
and since $x_d$ is superficial for $(\x)$, we must have
$s \in (x_1, \dots, x_d)^{m-1}$. Let us write
$$s = \sum_{|\n| = m-1} a_\n \x^\n.$$
We can now write (\ref{eq361}) as
\begin{align}
\sum_{|\n| = m-1} a_\n \x^\n x_d - \sum_{|\n| = m} d_\n \x^\n = 0, \label{eq362}
\end{align}
where $Q(x_1, \dots, x_{d-1}) = \sum_{|\n| = m} d_\n \x^\n \in (x_1, \dots, x_{d-1})^m$. Let $P(x_1, \dots, x_d) = \sum_{|\n| = m-1} a_\n \x^\n.$ Then, (\ref{eq362}) gives a relation on $x_1, \dots, x_d$ of degree $m$, namely
$$H(T_1, \dots, T_d) = T_d P(T_1, \dots, T_d) - Q(T_1, \dots, T_{d-1}).$$
Let $r_{\m_0} \T^{\m_0}$ (where $|\m_0| = m$) be the leading term of $H(T_1, \dots, T_d)$. It follows from Lemma \ref{map} that
$$H(x_1, \dots, x_d) - H(x_1, \dots, x_{d-1}, y_d) \in (\y^\m ~|~ \m < \m_0).$$
Moreover, $H(x_1, \dots, x_d) - H(x_1, \dots, x_{d-1}, y_d) = x_d P(x_1, \dots, x_d) - y_d P(x_1, \dots, x_{d-1}, y_d)$. Thus
\begin{align}
x_d P(x_1, \dots, x_d) - y_d P(x_1, \dots, x_{d-1}, y_d) \in (\y^\m ~|~ \m < \m_0). \label{eq363}
\end{align}
Since the calculation done in Lemma \ref{map} is formal on the coefficients of the relation $F(T_1, \dots, T_d)$, in our situation it is formal on the coefficients of $T_d P(T_1, \dots, T_d)$, which are exactly the same as those of $P(T_1, \dots, T_d)$. Therefore, the same calculation as in (\ref{eq363}) would hold for $P(x_1, \dots, x_d)$. Hence, for $\m_0' = \m_0 - (0, \dots, 0, 1)$, we have
\begin{align}
P(x_1, \dots, x_d) - P(x_1, \dots, x_{d-1}, y_d) \in (\y^\m ~|~ \m < \m_0'). \label{eq364}
\end{align}
It now follows from (\ref{eq364}) that
$$s = P(x_1, \dots, x_d) \in (x_1, \dots, x_{d-1}, y_d)^{m-1}.$$
The lemma is proved.
\end{proof}

\begin{lemma} \label{superficial}
Let $\x$ and $\y$ be as in Lemma \ref{super} in a local ring
$(R,\mm)$ with infinite residue field. There exist $y'_1, \dots,
y'_{d-1}$ such that $(y'_1, \dots, y'_{d-1}, y_d)=(\y)$, and $y_d,
y'_{d-1},\dots, y'_1$ is a superficial sequence.
\end{lemma}

\begin{proof}
Take $y'_1, \dots, y'_{d-1}$ to be general linear combinations of
$x_1,\dots, x_{d-1}$.
\end{proof}


\section{Ramsey numbers} \label{ramsey-lemma}

In this section, we provide a ``Ramsey number'' combinatorial lemma which will be used to establish uniform bounds on relation type of parameter ideals in the next two sections. 

For a set $S$, and a positive integer $l$, we denote by $[S]^l$ the set of all subsets of $l$ elements of $S$. We shall use the following infinite version of Ramsey's theorem \cite{ramsey} (see also \cite[Theorem 3.4]{adhikari}).

\begin{lemma} \label{ramsey}
Let $n$ and $l$ be two given positive integers. Let $T = \{ x_1, x_2, \dots \}$ be an infinite countable set. Then for any way of colouring $[T]^l$ using $n$ colours, there is an infinite subset $U$ of $T$ with all its subsets of $l$ elements having the same colour.
\end{lemma}

For any $d$-uple $A = (a_1, \dots, a_d)$ of nonnegative integers, we denote $|A| = \sum_{j=1}^d a_j$. For two $d$-uples of nonnegative integers $A = (a_1, \dots, a_d)$ and $B = (b_1, \dots, b_d)$, we write $A \preceq B$ if and only if $a_j \le b_j$ for all $j =1, \dots, d$. A {\sl chain} of length $l$ is a sequence of $d$-tuples of nonnegative integers $A_1 \preceq A_2 \preceq \dots \preceq A_l$. For an $d$-uple $A = (a_1, \dots, a_d)$ and a number $t \in \{ 1, \dots, d \}$, we use $A(t)$ to denote the $t$-th entry of $A$, i.e. $A(t) = a_t$.

\begin{lemma} \label{ramseynumber}
Suppose $d$ and $l$ are given positive integers, and $k$ is a nonnegative integer. Then, there exists a positive number $M = M(d,k,l)$ such that for any sequence of $d$-uples of nonnegative integers $\A = \langle A_1, A_2, \dots, A_M \rangle$, in which $|A_i| = k+i$, we can always find a chain $A_{i_1} \preceq A_{i_2} \preceq \dots \preceq A_{i_l}$ of length $l$ with $1 \le i_1 < i_2 < \dots < i_l \le M$.
\end{lemma}

\begin{proof} By contradiction, suppose the assertion is not true. That is, for any $M \ge k$, there is a sequence $\A_M = \langle A_{M,1}, A_{M,2}, \dots, A_{M,M} \rangle$ of $d$-uples of nonnegative integers such that $|A_{M,i}| = k+i$ for $i = 1, \dots, M$, which has no subchain of length $l$.

We shall first inductively construct an infinite sequence of infinite subsets $\M_1 \supseteq \M_2 \supseteq \dots \supseteq \M_n \supseteq \dots$ of $\NN$ as follows. Since $|A_{M,1}| = k+1$ for all $M \ge k$, there must be an infinite subset $\M_1$ of $\NN$ such that for any $M, N \in \M_1$, we have $A_{M,1} = A_{N,1}$. Suppose $\M_1, \dots, \M_i$ ($i \ge 1$) have been constructed. Since for each $M \in \M_i$, $|A_{M,i+1}| = k+i+1$, there must be an infinite subset $\M_{i+1}$ of $\M_i$ such that for any $M,N \in \M_{i+1}$, we have $A_{M,i+1} = A_{N,i+1}$.

Let $\A = \langle A_1, A_2, \dots \rangle$ be the sequence defined by letting $A_i = A_{M,i}$, where $M$ is an arbitrary element of $\M_i$. Clearly, by definition, $|A_i| = k+i$.

\begin{claim} \label{ramseyclaim}
The sequence $\A = \langle A_1, A_2, \dots \rangle$ does not have a subchain of length $l$.
\end{claim}

{\noindent \sc Proof of Claim.} Indeed, if $B_1 \preceq B_2 \preceq \dots \preceq B_l$ is a subchain of $\A$, and $B_l = A_s$, then for any $M \in \M_s$, the sequence $\A_M = \langle A_{M,1}, A_{M,2}, \dots \rangle$ contains the subchain $B_1 \preceq B_2 \preceq \dots \preceq B_l$ of length $l$. This is a contradiction. \qed

Now, let us consider a colouring of $[\A]^2$ using $(d+1)$ colours as follows. Suppose $i < j$. If $A_i \preceq A_j$, then we colour $\{ A_i, A_j \}$ by 0. Otherwise, there must be an integer $t \in \{ 1, \dots, n \}$ such that $A_i(s) \le A_j(s)$ for any $s < t$ and $A_i(t) > A_j(t)$, and we colour $\{ A_i, A_j \}$ by $t$. Clearly, this is a valid $(d+1)$-colouring of $[\A]^2$. By Ramsey's theorem (Theorem \ref{ramsey}), there exists an infinite subset $\U = \{ U_1, U_2, \dots \}$ of $\A$ such that $[\U]^2$ has one colour. Suppose $[\U]^2$ is coloured by $c$. If $c > 0$, then we obtain an infinite sequence of nonnegative integers $U_1(c) > U_2(c) > \dots $, which is impossible. Thus, $c = 0$, and so $U_1 \preceq U_2 \preceq \dots $ is a chain of infinite length in $\A$. This contradicts Claim \ref{ramseyclaim}.

Hence, we always get a contradiction. The lemma is proved.
\end{proof}


\section{Relation type in generalized Cohen-Macaulay rings} \label{bound-rt}

The goal of this section will be to give a new argument that rings of
finite local cohomology have uniformly bounded relation type. Throughout this section, $(R,\mm,k)$ will denote a local ring $R$ with maximal ideal $\mm$ and residue field $k$. 

We recall that a Noetherian local ring $R$ of dimension $d$ is
said to have {\sl finite local cohomology} (f.l.c.) if
$H_{\mm}^i(R)$ is finitely generated for $i=0,\dots,d-1$ (and hence is of finitely length). Rings
with finite local cohomology are called {\sl generalized
Cohen-Macaulay}. We observe that $R$ has f.l.c. if and only if $\hat R$ has f.l.c. if and only if $\dim \hat R/p=\dim\hat R$ for every minimal prime $p$
of $\hat R$, and $\hat R_p$ is Cohen-Macaulay for all $p\neq \hat
\mm$.

Notice that if a Notherian local ring $R$ has f.l.c. then $\A(R)$
is $\mm$-primary.

We shall start with a result in dimension 2 through which the argument in the general situation becomes more transparent. Notice that the bound $l_1+l_0$ of Theorem~\ref{2parameters} improves the bound $2l_1+l_0$ of \cite[Theorem 4.1]{Wang:1997}.

\begin{thm}\label{2parameters}
Let $(R,\mm,k)$ be an equidimensional local ring of dimension $2$ such that
$\length( H^0_\mm(R))=l_0 < \8$ and $\length( H^1_\mm(R))=l_1 < \infty$.
Then $R$ has a uniform bound $l_1+l_0$ on relation type of
parameter ideals.
\end{thm}

\begin{proof} Let $I = (x,y)$ be a parameter ideal of $R$. By
\cite[Lemma 2.2]{Wang:1997} we can assume that $H^0_\mm(R)=0$ and show
that $\rt(I) \leq l_1$. Suppose
$$F(T_1, T_2) =r_NT_1^N+r_{N-1}T_1^{N-1}T_2+\dots + r_0T_2^N,$$
with $N>l_1$, provides a relation on $(x,y)$. That is,
\begin{align}
r_Nx^N+r_{N-1}x^{N-1}y+\dots + r_0y^N=0. \label{2eq}
\end{align}
We may assume that $r_N \not= 0$, otherwise we can factor out a power of $T_2$.

Let $\gamma \in \A(R)$ be a homology multiplier which is part of a s.o.p. (in particular, $\gamma$ is a non-zero-divisor). Let $l=\length(H^0_\mm(R/\gamma R))$. By \cite[Lemma 3.7]{Wang:1997} we have
that $l\leq l_1.$
It follows from (\ref{2eq}) that $r_N \in y:x^N$. Thus, since $\gamma$ is a homology multiplier, we have
$$\gamma r_N \in (y).$$
Similarly, (\ref{2eq}) implies that
$r_Nx+r_{N-1}y \in y^2:x^{N-1}$, and so
$$\gamma (r_Nx+r_{N-1}y) \in
(y^2).$$
Proceeding in this way, we obtain a sequence of relations as follows.
\begin{eqnarray}
\left\{ \begin{array}{lcl} \gamma r_N & = & s_1y, \\
\gamma (r_Nx+r_{N-1}y) & = & s_2y^2, \\
& \vdots & \\
\gamma (r_Nx^{N-1}+r_{N-1}x^{N-2}y+\dots+r_1y^{N-1}) & = & s_Ny^N, \end{array}
\right. \label{2steps}
\end{eqnarray}
where $s_1,\dots,s_N \in R$.

For every $1\leq j\leq N$, from (\ref{2steps}) we have $s_j \in (\gamma:y^j) \subseteq H^0_\mm(R/\gamma R)$. Since $N > l_1$, there exists $p \le l_1$ such that
$$s_{p+1} \in (\gamma, s_1,\dots, s_p).$$ Write
$$s_{p+1}=a\gamma+b_1s_1+\dots+b_ps_p.$$
Substituting this into the relation
$\gamma(r_Nx^p+r_{N-1}x^{p-1}y+\dots+r_{N-p}y^p)=s_{p+1}y^{p+1}$
in (\ref{2steps}), we obtain
\begin{align*}
\gamma(r_Nx^p+r_{N-1}x^{p-1}y+\dots+r_{N-p}y^p) & = (a\gamma+b_1s_1+\dots+b_ps_p)y^{p+1} \nonumber \\
& = a\gamma y^{p+1} + b_1s_1y^{p+1} + \dots + b_ps_py^{p+1}.
\end{align*}
i.e.
\begin{align}
\gamma\Big(r_N x^p + r_{N-1}x^{p-1}y + \dots + (r_{N-p}-ay)y^p\Big) = (b_1y^p)s_1y + (b_2y^{p-1})s_2y^2 + \dots + (b_py)s_py^p. \label{2relation}
\end{align}
We can use (\ref{2steps}) to replace $(b_{j+1} y^{p-j})s_{j+1}y_{j+1}$ by $\gamma\Big((r_Nb_{j+1}y^{p-j})x^j + \dots + (r_{N-j}b_{j+1}y^{p-j})y^j\Big)$ in (\ref{2relation}), for $j=0, \dots, p-1$. Observe that for each $j = 1, \dots, p-1$, we have $(r_N b_{j+1} y^{p-j}) x^j = (r_Nb_{j+1}) x^jy^{p-j} < r_N x^p$. Thus, after moving everything to the left hand side and factoring out $\gamma$, (\ref{2relation}) gives us
$$\gamma\Big(r_Nx^p+ {\rm smaller \ terms}\Big)=0.$$
This implies, since $\gamma$ is a non-zero-divisor, that
$$ r_Nx^p + \ \text{smaller terms} \ = 0.$$
Therefore, we get a new relation on $(x,y)$,
$$G(T_1, T_2) = r_NT_1^p+ {\rm smaller \ terms}$$
which has the same leading coefficient as $F(T_1, T_2)$ but is of lower degree.

Now, write $F = T_1^{N-p}G(T_1, T_2)+T_2H(T_1, T_2)$. Then, clearly $H(T_1, T_2)$ also provides a relation on $(x,y)$ and is of smaller degree than $F(T_1, T_2)$. That is, $F(T_1, T_2)$ can be written as a combination of relations of lower degrees. This proves our result.
\end{proof}

To prove the result for generalized CM rings of any dimension, we shall need the following lemma of Schenzel (cf. \cite[Theorem 3.2]{Wang:1997}).

\begin{lemma} \label{gCMcohombound}
Let $(R,\mm,k)$ be a $d$-dimensional Noetherian ring having finite local cohomology. Let $(x_1, \dots, x_s)$ be part of a s.o.p. in $R$. Then,
$$\lambda \Big( \sfrac{(x_2, \dots, x_s):x_1}{(x_2, \dots, x_s)} \Big) \le \sum_{i=0}^{s-1} {s-1 \choose i} \lambda \big( H^i_\mm(R) \big).$$
\end{lemma}

\begin{remark} \label{gCMh0}
Under the assumptions of Lemma \ref{gCMcohombound}, we have
$$H^0_\mm \big( R/(x_2, \dots, x_s) \big) = \sfrac{(x_2, \dots, x_s):x_1^{\infty}}{(x_2, \dots, x_s)} = \sfrac{(x_2, \dots, x_s):x_1^n}{(x_2, \dots, x_s)},$$
for some positive integer $n$.
Since if $(x_1, \dots, x_s)$ is part of a s.o.p. then so is $(x_1^n, x_2, \dots, x_s)$, Lemma \ref{gCMcohombound} now gives
$$\lambda \big( H^0_\mm \big(R/(x_2, \dots, x_s) \big) \big) \le \sum_{i=0}^{s-1} {s-1 \choose i} \lambda \big( H^i_\mm(R) \big).$$
\end{remark}

Our result for generalized CM rings of any dimension is stated as follows.

\begin{thm} \label{dparameters}
Let $(R,\mm,k)$ be a local ring of dimension $d$ such that
$\length( H^i_\mm(R)) < \8$ for $0\leq i\leq d-1$. Then $R$ has a
uniform bound on relation type of parameter ideals.
\end{thm}

\begin{proof} By \cite[Lemma 4.1]{Wang1}, since our hypotheses pass to the completion, we may assume that $R$ is complete. By \cite[Lemma 2.2]{Wang:1997}, we may also assume that $H^0_\mm(R) = 0$. Furthermore, by passing to a faithful extension of $R$ if necessary, we may assume that $k$ is infinite (see the proof of \cite[Theorem 4.4]{Wang:1997}).

Let $I=(x_1,\dots,x_d)$ be a parameter ideal in $R$. We can pick $x_d, x_{d-1}, \dots, x_1$ to form a superficial sequence and $x_d$ not to be a zero-divisor. We will use graded reverse lex monomial ordering with $x_1 > x_2 > \dots > x_d$ and $T_1 > T_2 > \dots > T_d$. For $q = 2, \dots, d$, let $L_q = \big[\sum_{i=0}^{d-q+1} {d-q+1 \choose i} \lambda \big( H^i_\mm(R) \big)\big]+1$. Let $L = (L_2-1) + \dots + (L_d-1) + 1$. We shall recursively construct the following sequence of finite numbers: let $M(d,k,l)$ denote the Ramsey number determined as in Lemma \ref{ramseynumber}, let $K_1 = M(d,1,L) = N_1$; and for $i \ge 2$, let $M_i = \sum_{l=0}^{K_{i-1}} {l+i-2 \choose i-2}$ (observe that $M_i$ is the number of all $(i-1)$-tuples of non-negative integers whose sum is at most $K_{i-1}$; we shall need this fact later on in the proof), $N_i = M(d, 2K_{i-1}, M_i(L-1)+1)$ and $K_i = 2K_{i-1}+N_i$.

To get the conclusion, it suffices to prove that any relation on $x_1, \dots, x_d$ of degree greater than $K_{d-1}$ can be written as a combination of relations of smaller degrees. Consider an arbitrary relation on degree $N > K_{d-1}$ in $x_1, \dots, x_d$
$$F = F(T_1, \dots, T_d) = \sum_{|\n| = N} r_\n \T^\n.$$
That is,
\begin{eqnarray}
F(x_1, \dots, x_d) = \sum_{|\n| = N} r_\n \x^\n = 0. \label{startrel}
\end{eqnarray}
We can assume that $T_d$ does not divide the leading term of $F$, since otherwise we can factor out $T_d$ and get a relation on smaller degree (since $x_d$ is not a zero-divisor). Let $\ord_i(\T^\n)$ be the $i$-th component of $\n$, and let $\ord_{< j}(\T^\n) = \sum_{i < j} \ord_i(\T^\n)$. We shall show that for any $1 \le j \le d-1$, if $F$ contains a term $r_\m \T^\m$ (say $\m = (m_1, \dots, m_d)$) with $\ord_{< j+1}(\T^\m) > K_j$ then we can write $F(T_1, \dots, T_d) = H(T_1, \dots, T_d)G(T_1, \dots, T_d) + F'(T_1, \dots, T_d)$, where $H$ is a monomial divisible by $\prod_{i > j} T_i^{m_i}$, $G$ and $F'$ both provide relations in $x_1, \dots, x_d$, the leading term $r_\k \T^\k$ of $G$ satisfies the condition $r_\k \T^\k \big| r_\m \T^\m, \ord_{< j+1} (\T^\k) \le K_j$ and $\ord_{j+1}(\T^\k) = \dots = \ord_d(\T^\k) = 0$ (in particular, the degree of $G$ is bounded by $K_j$), and all terms $r_{\n_1}\T^{\n_1}$ of $F'$ with $\ord_{< j+1}(\T^{\n_1}) > K_j$ are smaller than $r_\m \T^\m$. By taking $j = d-1$, and successively eliminate terms $r_\n \T^\n$ of $F$ with $\ord_{< d}(\T^\n) > K_{d-1}$, we then prove our theorem.

Suppose our assertion is not true. Let $j$ be the smallest index for which our assertion fails. Let $F$ be a relation such that our assertion fails for this value of $j$. In particular, $F$ contains a term $r_\n \T^\n$ such that $\ord_{< j+1}(\T^\n) > K_j$. Let $r_{\n_0} \T^{\n_0}$ be the largest term of $F$ for which $\ord_{< j+1}(\T^{\n_0}) > K_j$. Among all such relations $F$'s for which our assertion fails for $j$, we shall pick $F$ such that $r_{\n_0} \T^{\n_0}$ is smallest possible. Let $\n_0 = (n_{01}, \dots, n_{0d})$. For simplicity, we write $r$ for $r_{\n_0}$. We shall derive a contradiction.

From the choice of $j$, we may assume that $K = \sum_{i=1}^{j-1} n_{0i} \le K_{j-1}$. We first observe that
\begin{align}
n_{0j} > K_j - K_{j-1} = K_{j-1}+N_j. \label{gCMn0}
\end{align}
Let $\J$ denote the set of all monomials in $T_1, \dots, T_d$ that appear in the expression of $F(T_1, \dots, T_d)$. Let $P(T_1, \dots, T_d)$ be the sum of all terms of $F$ that are divisible by $\prod_{i > j}T_i^{n_{0i}}$, i.e.
$$ P(T_1, \dots, T_d) = \prod_{i > j}T_i^{n_{0i}} \Big(r \sfrac{\T^{\n_0}}{\prod_{i > j}T_i^{n_{0i}}} + \sum_{\n \in \J, \n \not= \n_0, ~ \prod_{i > j}T_i^{n_{0i}} \big| \T^\n} r_{\n} \sfrac{\T^\n}{\prod_{i > j}T_i^{n_{0i}}} \Big). $$
We first observe that if $\n > \n_0$, then from the choice of $r_{\n_0} \T^{\n_0}$, we must have
$\ord_{< j+1}(\T^\n) \le K_j = \ord_{< j+1}(\T^{\n_0}),$
whence
\begin{align}
\sum_{i > j} \ord_i(\T^\n) > \sum_{i > j} \ord_i(\T^{\n_0}). \label{gCMord}
\end{align}
This implies that for all $\T^\n \in \J$ such that $\n > \n_0$, we must have $\prod_{i > j}T_i^{n_{0j}} \nmid ~ \T^\n$, i.e. $r_\n \T^\n$ is not in $P(T_1, \dots, T_d)$.
It now follows from (\ref{startrel}) that
\begin{align}
r \sfrac{\x^{\n_0}}{\prod_{i > j}x_i^{n_{0i}}} + \sum_{\n \in \J, \n < \n_0, ~\prod_{i > j}x_i^{n_{0i}} \big| \x^\n} r_{\n} \sfrac{\x^\n}{\prod_{i > j}x_i^{n_{0i}}} & = \sfrac{P(x_1, \dots, x_d)}{\prod_{i > j}x_i^{n_{0i}}} \nonumber \\
& \in \Big(\sum_{\n \in \J, ~\prod_{i > j}x_i^{n_{0i}} \nmid ~ \x^\n} r_\n \x^\n\Big) : \prod_{i > j}x_i^{n_{0i}}. \label{gCMcolon}
\end{align}
Since $R$ is generalized Cohen-Macaulay, i.e. $\A(R)$ is $\mm$-primary, there exists a positive integer $q_j'$ such that $x_j^{q_j'} \in \A(R)$. Let $B = (x_{j+1}, \dots, x_d)$. It can be seen that if $\n < \n_0$ and $\prod_{i > j}T_i^{n_{0j}} \nmid ~ \T^\n$ then there exists $l > j$ such that $\ord_l(\T^\n) > \ord_l(\T^{\n_0})$. On the other hand, if $\n > \n_0$ then it follows from (\ref{gCMord}) that there also exists $l > j$ such that $\ord_l(\T^\n) > \ord_l(\T^{\n_0})$. Thus, from (\ref{gCMcolon}), we have
\begin{align}
x_j^{q_j'} \Big(r \sfrac{\x^{\n_0}}{\prod_{i > j}x_i^{n_{0i}}} + \sum_{\n \in \J, \n < \n_0, ~\prod_{i > j}x_i^{n_{0i}} \big| \x^\n} r_{\n} \sfrac{\x^\n}{\prod_{i > j}x_i^{n_{0i}}}\Big) \in B. \label{red1}
\end{align}
Since $x_d, \ldots, x_1$ form a superficial sequence, (\ref{red1})
and Remark \ref{sup-high-deg} imply there exists an integer
$q_j$ such that
\begin{align}
x_j^{q_j} \Big(r \sfrac{\x^{\n_0}}{\prod_{i > j}x_i^{n_{0i}}} + \sum_{\n \in \J, \n < \n_0, ~\prod_{i > j}x_i^{n_{0i}} \big| \x^\n} r_{\n} \sfrac{\x^\n}{\prod_{i > j}x_i^{n_{0i}}}\Big) \in B(x_1, \dots x_d)^{q_j+K+n_{0j}-1}. \label{red2}
\end{align}
Let us write (\ref{red2}) as
\begin{align}
r x_1^{n_{01}} \dots x_{j-1}^{n_{0(j-1)}} x_j^{n_{0j}+q_j} + \sum_{|\m| = q_j+K+n_{0j}, \m < (n_{01}, \dots, n_{0j}+q_j, 0, \dots, 0)} u_\m \x^\m = 0, \label{neweq}
\end{align}
and let
\begin{align}
Q(T_1, \dots, T_d) & = r T_1^{n_{01}} \dots T_{j-1}^{n_{0(j-1)}} T_j^{n_{0j}+q_j} + \sum_{|\m| = q_j+K+n_{0j}, \m < (n_{01}, \dots, n_{0j}+q_j, 0, \dots, 0)} u_\m \T^\m \nonumber \\
& = r \T^{\m_0} + \sum_{|\m| = q_j+K+n_{0j}, \m < \m_0} u_\m \T^\m, \label{neweq0}
\end{align}
where $\m_0 = (m_{01}, \dots, m_{0d}) = (n_{01}, \dots, n_{0(j-1)}, n_{0j}+q_j, 0, \dots, 0)$. Then, $Q(T_1, \dots, T_d)$ gives a new relation on $x_1, \dots, x_d$.

Let $\a$ be the collection of all terms of $Q$ that are smaller than $r \T^{\m_0}$. For a term $u_\m \T^\m \in \a$, we observe the following. If $\ord_i(\T^\m) = 0$ for all $i > j$ then we must have $\ord_j(\T^\m) > \ord_i(\T^{\m_0}) = m_{0j}$, and so $\ord_{< j}(\T^\m) \le K$. Otherwise, suppose there exists $l > j$ such that $\ord_l(\T^\m) > 0$. If $\ord_{< j}(\T^\m) > K_{j-1}$ then from the choice of $j$, we may write $Q(T_1, \dots, T_d) = H(T_1, \dots, T_d)G(T_1, \dots, T_d) + Q'(T_1, \dots, T_d)$ where $H$ is a monomial divisible by $T_l$ (since $\ord_l(\T^\m) > 0$), $G$ and $Q'$ are relations on $x_1, \dots, x_d$, $G$ has the form
\begin{align}
G(T_1, \dots, T_d) = v_{\k_0} \T^{\k_0} + \sum_{\k < \k_0, |\k| = |\k_0|} v_\k \T^\k, \label{simeq}
\end{align}
with $v_{\k_0} \T^{\k_0} \big| u_\m \T^\m$, $\ord_{< j}(\T^{\k_0}) \le K_{j-1}$ and $\ord_j(\T^{\k_0}) = \dots = \ord_d(\T^{\k_0}) = 0$. Notice that since $H$ is divisible by $T_l$, all terms in $H(T_1, \dots, T_d)G(T_1, \dots, T_d)$ are smaller than $r \T^{\m_0}$. By replacing $Q$ by $Q'$ and repeating the process to successively remove all terms $u_\m \T^\m \in \a$ such that $\ord_{< j}(\T^\m) > K_{j-1}$ from $Q$, we may assume that in our relation $Q(T_1, \dots, T_d)$, every term $u_\m \T^\m$ satisfies the condition $\ord_{< j}(\T^\m) \le K_{j-1}$.

Let $\gamma \in \A(R)$ be a homology multiplier such that $(\gamma, x_1, \dots, \widehat{x_j}, \dots, x_d)$ is a s.o.p. Since $\gamma$ is part of a s.o.p. and $H^0_\mm(R) = 0$, $\gamma$ is a non-zero-divisor. Let $\MM$ denote the set of all monomials appearing in $Q$. Recall that $K = \sum_{i < j} n_{0i} = \sum_{i < j} m_{0i} \le K_{j-1}$ and $B = (x_{j+1}, \dots, x_d)$.

Let $F_1(T_1, \dots, T_d)$ be the sum of all terms of $Q$ that are divisible by $T_j^{m_{0j}-K_{j-1}}$, i.e.
\begin{align*}
F_1(T_1, \dots, T_d) = T_j^{m_{0j}-K_{j-1}} \Big( r \sfrac{\T^{\m_0}}{T_j^{m_{0j}-K_{j-1}}}
+ \sum_{\m \in \MM, \m < \m_0, ~ T_j^{m_{0j}-K_{j-1}} \big| \T^\m} u_{\m} \sfrac{\T^\m}{T_j^{m_{0j}-K_{j-1}}} \Big).
\end{align*}
Let $G_1(T_1, \dots, T_d) = \sfrac{F_1(T_1, \dots, T_d)}{T_j^{m_{0j}-K_{j-1}}}$. Then, (\ref{neweq}) gives
\begin{align*}
G_1(x_1, \dots, x_d) & \in \Big(\x^\m \in \MM ~\Big|~ \m < \m_0, \ord_j(\x^\m) < m_{0j}-K_{j-1} \Big) : x_j^{m_{0j}-K_{j-1}} \\
& = \Big( \x^\m \in \MM ~\Big|~ \ord_j(\x^\m) < m_{0j}-K_{j-1}, \ \text{and} \\
& \hspace{15ex} \sum_{i > j} \ord_i(\x^\m) > K \Big) : x_j^{m_{0j}-K_{j-1}}.
\end{align*}
This is because $\sum_{i > j} \ord_i(\x^\m) > \sum_{i=1}^d \ord_i(\x^\m) - \ord_{< j}(\x^\m) - (m_{0j}-K_{j-1}) = \sum_{i=1}^d \ord_i(\x^{\m_0}) - \ord_{< j}(\x^\m) - (m_{0j}-K_{j-1}) \ge K+m_{0j} - K_{j-1} - (m_{0j}-K_{j-1}) = K$. Therefore, since $\gamma$ is a homology multiplier, we have
$$\gamma G_1(x_1, \dots, x_d) \in (x_1, \dots, \widehat{x_j}, \dots, x_d)^{K_{j-1}}B^{K+1}.$$
Similarly, let $F_2(T_1, \dots, T_d)$ be the sum of all terms of $Q$ that are divisible by $T_j^{m_{0j}-K_{j-1}-1}$, i.e.
\begin{align*}
F_2(T_1, \dots, T_d) = T_j^{m_{0j}-K_{j-1}-1} \Big( r \sfrac{\T^{\m_0}}{T_j^{m_{0j}-K_{j-1}-1}} + \sum_{\m \in \MM, \m < \m_0, ~ T_j^{m_{0j}-K_{j-1}-1} \big| \T^\m} u_{\m} \sfrac{\T^\m}{T_j^{m_{0j}-K_{j-1}-1}} \Big),
\end{align*}
and let $G_2(T_1, \dots, T_d) = \sfrac{F_2(T_1, \dots, T_d)}{T_j^{m_{0j}-K_{j-1}-1}}$, we then have
\begin{align*}
G_2(x_1, \dots, x_d) & \in \Big( \x^\m \in \MM ~\Big|~ \m < \m_0, \ord_j(\x^\m) < m_{0j}-K_{j-1}-1 \Big) : x_j^{m_{0j}-K_{j-1}-1} \\
& = \Big( \x^\m \in \MM ~\Big|~ \ord_j(\x^\m) < m_{0j}-K_{j-1}-1 \ \text{and} \\
& \hspace{15ex} \sum_{i > j} \ord_i(\x^\m) > K+1 \Big) : x_j^{m_{0j}-K_{j-1}-1}.
\end{align*}
This is because $\sum_{i > j} \ord_i(\x^\m) > \sum_{i=1}^d \ord_i(\x^\m) - \ord_{< j}(\x^\m) - (m_{0j}-K_{j-1}-1) = \sum_{i=1}^d \ord_i(\x^{\m_0}) - \ord_{< j}(\x^\m) - (m_{0j}-K_{j-1}-1) \ge K+m_{0j} - K_{j-1} - (m_{0j}-K_{j-1}-1) = K+1$. Thus, again since $\gamma$ is a homology multiplier, we get
$$ \gamma G_2(x_1, \dots, x_d) \in (x_1, \dots, \widehat{x_j}, \dots, x_d)^{K_{j-1}}B^{K+2}.$$
Proceed in this way, we obtain the following:
\begin{eqnarray}
\left\{ \begin{array}{rcl}
\gamma G_1(x_1, \dots, x_d) & \in & (x_1, \dots, \widehat{x_j}, \dots, x_d)^{K_{j-1}}B^{K+1}, \\
\gamma G_2(x_1, \dots, x_d) & \in & (x_1, \dots, \widehat{x_j}, \dots, x_d)^{K_{j-1}}B^{K+2}, \\
& \vdots & \\
\gamma G_p(x_1, \dots, x_d) & \in & (x_1, \dots, \widehat{x_j}, \dots, x_d)^{K_{j-1}}B^{K+p}, \\
& \vdots & \\
\gamma G_{m_{0j}-K_{j-1}+1}(x_1, \dots, x_d) & \in & (x_1, \dots, \widehat{x_j}, \dots, x_d)^{K_{j-1}}B^{K+m_{0j}-K_{j-1}+1},
\end{array} \right. \label{gCMincl}
\end{eqnarray}
where
$$G_p(T_1, \dots, T_d) = \sfrac{F_p(T_1, \dots, T_d)}{T_j^{m_{0j}-K_{j-1}-p+1}} = r T_1^{m_{01}} \dots T_{j-1}^{m_{0(j-1)}} T_j^{K_{j-1}+p-1} + \ \text{smaller terms}, \ $$
and $F_p(T_1, \dots, T_d)$ is the sum of all terms in $Q$ that are divisible by $T_j^{m_{0j}-K_{j-1}-p+1}$, for $1 \le p \le m_{0j}-K_{j-1}+1$.

Observe that $G_p(x_1, \dots, x_d)$
is of degree $(K_{j-1}+K+p-1)$ on $x_1, \dots, \widehat{x_j}, \dots, x_d$. Let $C_p = (x_1, \dots, \widehat{x_j}, \dots, x_d)^{K_{j-1}}B^{K+p}$. It follows from (\ref{gCMincl}) that for each $1 \le p \le m_{0j}-K_{j-1}+1$, we can write $\gamma G_p(x_1, \dots, x_d) = H_p(x_1, \dots, x_d)$, where $H_p(x_1, \dots, x_d) \in C_p$ is a polynomial on $x_1, \dots, \widehat{x_j}, \dots, x_d$ of degree $(K_{j-1}+K+p)$. From now on, we shall write $G_p$ and $H_p$ for $G_p(x_1, \dots, x_d)$ and $H_p(x_1, \dots, x_d)$, respectively. We order the terms in $H_p$ with respect to our monomial ordering (graded reverse lex), and let $s_p$ be the leading coefficient of $H_p$. That is,
\begin{eqnarray}
\left\{ \begin{array}{rcl}
\gamma G_1 & = & s_1 \x^{\m_1} + \sum_{\m < \m_1} s_\m \x^\m, \\
\gamma G_2 & = & s_2 \x^{\m_2} + \sum_{\m < \m_2} s_\m \x^\m, \\
& \vdots & \\
\gamma G_p & = & s_p \x^{\m_p} + \sum_{\m < \m_p} s_\m \x^\m, \\
& \vdots & \\
\gamma G_{m_{0j}-K_{j-1}+1} & = & s_{m_{0j}-K_{j-1}+1} \x^{\m_{m_{0j}-K_{j-1}+1}} + \sum_{\m < \m_{m_{0j}-K_{j-1}+1}} s_\m \x^\m,
\end{array} \right. \label{gCMsteps}
\end{eqnarray}
where $\x^{\m_p} \in C_p, |\m_p| = K_{j-1}+K+p$. In the $p$-th equality of  (\ref{gCMsteps}), among all different ways of writing $\gamma G_p = H_p$ where $G_p$ is a polynomial expression in $x_1, \dots, x_d$ of degree $(K_{j-1}+K+p-1)$ and leading term $r x_1^{m_{01}} \dots x_{j-1}^{m_{0(j-1)}} x_j^{K_{j-1}+p-1}$ and $H_p$ is a polynomial expression in $x_1, \dots, x_d$ of degree $(K_{j-1}+K+p)$, we shall choose the one with smallest possible leading term $s_p \x^{\m_p}$ on the right hand side.

\begin{claim} \label{gCMclaim}
There exists an integer $p \le N_j$ such that $H_p = 0$.
\end{claim}

\noindent {\sc Proof of Claim.} From (\ref{gCMn0}) we have $m_{0j}-K_{j-1}+1 \ge n_{0j}-K_{j-1}+1 \ge N_j+1$. By contradiction, suppose the assertion is false. That is, $s_1, \dots, s_{N_j}$ are all non-zero. Fix an integer $1 \le p \le m_{0j}-K_{j-1}+1$, and suppose the $x_i$'s that appear in $\x^{\m_p}$ are in $\{ x_i ~|~ i \ge h\} \backslash \{x_j\}$ (and $h$ is chosen to be the largest integer with this property). Then
\begin{align*}
s_p & \in \Big( \gamma, \{ \x^\m \in C_p ~ \big|~ \m < \m_p \} \Big) : \x^{\m_p} \\
& = \Big( \gamma, \{ \x^\m \in C_p ~|~ \ord_i(\x^\m) > \ord_i(\x^{\m_p}) \ \text{for some} \ i \ge h+1 \} \Big) : \x^{\m_p}.
\end{align*}
Choose $\alpha \in \A(R)$ such that $\big(\gamma, \{ x_{h+1}, \dots, x_d \} \backslash \{ x_j \}, \alpha\big)$ is part of a s.o.p. (since $x_h$ divides $\x^{\m_p}$, we have $h \not= j$, so this choice is possible). Then,
$$ \alpha s_p \in \big(\gamma, \{ x_{h+1}, \dots, x_d\} \backslash \{x_j\}\big), $$
i.e.
$$ s_p \in \big(\gamma, \{ x_{h+1}, \dots, x_d\} \backslash \{x_j\} \big) : \alpha.$$
Since $R$ is a generalized Cohen-Macaulay ring and $\big(\gamma, \{ x_{h+1}, \dots, x_d \} \backslash \{ x_j \} , \alpha\big)$ is part of a s.o.p., $\alpha$ avoids all associated primes of $\big(\gamma, \{ x_{h+1}, \dots, x_d \} \backslash \{ x_j \}\big)$ except possibly $\mm$. Hence, we have
$$ \overline{s_p} \in H^0_\mm\big(R\big/\big(\gamma, \{ x_{h+1}, \dots, x_d\} \backslash \{x_j\} \big)\big), $$
where $\overline{s_p}$ denotes the image of $s_p$ in $R\big/\big(\gamma, \{ x_{h+1}, \dots, x_d\} \backslash \{x_j\} \big)$.

Consider the sequence $\{ \m_1, \m_2, \dots, \m_{N_j} \}$. Since $N_j = M(d,2K_{j-1},M_j(L-1)+1) \ge M(d,K_{j-1}+K,M_j(L-1)+1)$ and $|\m_1| = K_{j-1}+K+1$, it follows from Lemma \ref{ramseynumber} that there exists an increasing subsequence of length $V = M_j(L-1)+1$, $\m_{l_1} \preceq \m_{l_2} \preceq \dots \preceq \m_{l_V}$ with $1 \le l_1 < l_2 < \dots < l_V \le N_j \le m_{0j}-K_{j-1}$. Since $C_p = (x_1, \dots, \widehat{x_j}, \dots, x_d)^{K_{j-1}}B^{K+p}$, we have $\ord_{< j}(\x^{\m_p}) \le K_{j-1}$ for any $1 \le p \le m_{0j}-K_{j-1}+1$. Thus, since $M_j$ is the number of $(j-1)$-tuples of non-negative integers whose sum is at most $K_{j-1}$, we can choose from the sequence $\m_{l_1} \preceq \m_{l_2} \preceq \dots \preceq \m_{l_V}$ a subsequence of length $L$, $\m_{i_1} \preceq \m_{i_2} \preceq \dots \preceq \m_{i_L}$ such that they have the same $j-1$ first components, i.e., for any $1 \le a, b \le L$ we have $(\ord_1(\x^{\m_{i_a}}), \dots, \ord_{j-1}(\x^{\m_{i_a}})) = (\ord_1(\x^{\m_{i_b}}), \dots, \ord_{j-1}(\x^{\m_{i_b}}))$. Since $L = (L_2-1) + \dots + (L_d-1) + 1$, we can choose from this sequence yet another subsequence $\m_{p_1} \preceq \m_{p_2} \preceq \dots \preceq \m_{p_{L_q}}$, for some $2 \le q \le d$, such that
$$\overline{s_{p_t}} \in H^0_\mm\big(R\big/\big(\gamma, \{ x_{q+1}, \dots, x_d \} \backslash \{x_j\} \big)\big),$$
for any $t = 1, 2, \dots, L_q$. It follows from Remark \ref{gCMh0} that
$$\lambda \big( H^0_\mm\big(R\big/\big(\gamma, \{ x_{q+1}, \dots, x_d \} \backslash \{x_j\} \big)\big) \big) < L_q, $$
and therefore, we can write
$$ s_{p_{L_q}} = \sum_{i > q, i \not= j} a_i x_i + a_j \gamma + \sum_{i=1}^{L_q-1}b_is_{p_i}. $$
Substitute this into the equality $\gamma G_{p_{L_q}} = H_{p_{L_q}}$ to get
\begin{align}
\gamma (G_{p_{L_q}} - a_j\x^{\m_{p_{L_q}}}) & = \sum_{i > q, i \not= j} a_ix_i \x^{\m_{p_{L_q}}} + \sum_{i=1}^{L_q-1}b_is_{p_i}\x^{\m_{p_{L_q}}} + \sum_{\m < \m_{p_{L_q}}} s_\m \x^\m \nonumber \\
& = \sum_{i > q, i \not= j} (a_ix_q) \sfrac{x^{\m_{p_{L_q}}}x_i}{x_q} + \sum_{i=1}^{L_q-1} \Big(b_i\sfrac{\x^{\m_{p_{L_q}}}}{\x^{\m_{p_i}}}\Big)s_{p_i}\x^{\m_{p_i}} + \sum_{\m < \m_{p_{L_q}}} s_\m \x^\m. \label{simplifying}
\end{align}
Since $\x^{\m_{p_{L_q}}} \in C_{p_{L_q}}$, (after absorbing extra powers into the coefficient if necessary) $a_j\x^{\m_{p_{L_q}}} < \ r x_1^{m_{01}} \dots x_{j-1}^{m_{0(j-1)}} x_j^{K_{j-1}+p_{L_q}-1}$. Also, for $i > q$,  $\sfrac{\x^{\m_{p_{L_q}}}x_i}{x_q} < \x^{\m_{p_{L_q}}}$. Furthermore, since $\m_{p_{L_q}}$ and $\m_{p_i}$ have the same $j-1$ first components, we have $\sfrac{\x^{\m_{p_{L_q}}}}{\x^{\m_{p_i}}} \in (x_{j+1}, \dots, x_d)$, which then implies that all terms of $(b_i\sfrac{\x^{\m_{p_{L_q}}}}{\x^{\m_{p_i}}}\Big)G_{p_i}$ are smaller than $r x_1^{m_{01}} \dots x_{j-1}^{m_{0(j-1)}} x_j^{K_{j-1}+p_{L_q}-1}$. Thus, after using (\ref{gCMsteps}) to replace $s_{p_i}\x^{\m_{p_i}}$ by $\big(\gamma G_{p_i} - \sum_{|\m| = |\m_{p_i}|, \m < \m_{p_i}} s_\m x^\m\big)$ in (\ref{simplifying}) and moving terms involving $\gamma$ to the left hand side, we get a relation
$$\gamma \Big(G_{p_{L_q}} + \ \text{terms smaller than} \ r x_1^{m_{01}} \dots x_{j-1}^{m_{0(j-1)}} x_j^{K_{j-1}+p_{L_q}-1}\Big) = H_{p_{L_q}}'.$$
That is,
\begin{eqnarray}
\gamma \Big(r x_1^{m_{01}} \dots x_{j-1}^{m_{0(j-1)}} x_j^{K_{j-1}+p_{L_q}-1}+ \ \text{smaller terms}\Big) = H_{p_{L_q}}', \label{newrel}
\end{eqnarray}
where the leading term of $H_{p_{L_q}}'$ is strictly smaller than $s_{p_{L_q}}\x^{\m_{p_{L_q}}}$, which is the leading term of $H_{p_{L_q}}$. This contradicts the way (\ref{gCMsteps}) was chosen. Our claim is proved. \qed

We have just shown that there exists an integer $p \le N_j$ such that the $p$-th equality in (\ref{gCMsteps}) is
$ \gamma G_p = 0.$
Thus, since $\gamma$ is a non-zero-divisor, $G_p = 0$. That is,
$$ r x_1^{m_{01}} \dots x_{j-1}^{m_{0(j-1)}} x_j^{K_{j-1}+p-1}+ \ \text{smaller terms} \ = 0. $$
This gives a relation on $x_1, \dots, x_d$, which by abusing language we shall denote by $G_p$,
\begin{align*}
G_p(T_1, \dots, T_d) & = r T_1^{m_{01}} \dots T_{j-1}^{m_{0(j-1)}} T_j^{K_{j-1}+p-1} + \ \text{smaller terms} \\
& = r T_1^{n_{01}} \dots T_{j-1}^{n_{0(j-1)}} T_j^{K_{j-1}+p-1} + \ \text{smaller terms}.
\end{align*}
Observe that $K_{j-1}+p-1 \le K_{j-1}+N_j-1 < n_{0j}$ by (\ref{gCMn0}). We now can write
$$F(T_1, \dots, T_d) = \big[T_j^{n_{0j}-K_{j-1}-p+1} \prod_{i > j} T_i^{n_{0i}}\big] G_p(T_1, \dots, T_d) + F'(T_1, \dots, T_d),$$
where both $G_p$ and $F'$ are relations on $x_1, \dots, x_d$. It is clear that every term $r_\n \T^\n$ of $F'$ with $\ord_{< j+1}(\T^\n) > K_j$ is smaller than $r_{\n_0}\T^{\n_0}$. We obtain a contradiction. Hence, the theorem is proved.
\end{proof}


\section{Relation type in rings with non-Cohen-Macaulay locus of dimension one} \label{sectionCM1}

This section is devoted to treating the unknown situation where the ring $R$ has non-CM locus of dimension 1. Our main theorem shows that $R$ satisfies bounded relation type provided $\A(\hat{R})$ is a prime ideal in the completion $\hat{R}$ of $R$. As before, throughout the section, $(R,\mm,k)$ shall denote a local ring $R$ with maximal ideal $\mm$ and residue field $k$. 

\begin{lemma}\label{CM1cohombound}
Let $(R,\mm,k)$ be a complete unmixed local ring of dimension $d\ge
3$.  Assume that  $\A(R) = P$ is a dimension one prime. Then
$\lambda(H^1_\mm(R)) < \infty.$ Let $x_d \in R$ be a
non-zero-divisor such that its image in $R/P$ is in $\mm R/P -
\mm^\r R/P$ for some positive integer $\r$. Then for any system of
parameters $x_1,\ldots, x_d$ and any integer $1<j\leq d$ there is
a bound on $\length\left(H^0_\mm\left(R/(x_j,\ldots,
x_d)\right)\right)$, depending only on $\r$ and $j$.
\end{lemma}

\begin{proof}  The assertion that $\lambda(H^1_\mm(R)) < \infty$ follows
by duality.  Let $S \subseteq R$ be a Gorenstein ring  with $R$ module-finite
over $S$.  Then $H^1_\mm(R)$ is dual to $\Ext_S^{d-1}(R,S)$.  For any
non-maximal prime $Q \subseteq R$ of height $h$,  let $q = Q \cap S$. Then
$(\Ext_S^{d-1}(R,S))_Q = \Ext_{S_q}^{d-1}(R_Q, S_q)$ is dual to
$H_{QR_Q}^{h-d+1}(R_Q)$.  If $h < d-1$ this clearly vanishes.  Otherwise
$h = d-1$, and the assumption that $R$ is unmixed gives
$H_{QR_Q}^{0}(R_Q) = 0$.  Thus $\lambda(H^1_\mm(R)) = \lambda(\Ext_S^{d-1}(R,S))
<\infty$.

For simplicity of notation let $x_d=x$. We next observe that the
hypotheses give the non-CM locus of $R$ is $\{P, \mm\}$.  Since $x
\notin P$, the ring $R/xR$ is generalized CM. Thus, if we can bound
the lengths of $H^i_\mm(R/xR)$ for $0 \le i \le d-2$, depending
only on $i$ and $\r$ then Schenzel's result (Lemma \ref{gCMcohombound})
may be applied.

Since $P = \A(R)$ annihilates $H^i_\mm(R)$ for $0 \le i \le d-1$,
we have that $H^i_\mm(R)$ is an Artinian $R/P$ module.  Let $E =
E_{R/P}(R/\mm)$ be an injective hull. For each $1 < i < d$ there
is an exact sequence $0 \to H^i_\mm(R) \to E^{t_i} \to C_i \to 0$
(and $C_i$ is Artinian).  From the snake lemma applied to
multiplication of this short exact sequence by $x$ we get an exact
sequence
\begin{align}
0 \to \Ann_{H^{i}_\mm(R)} x \to (\Ann_E x)^{t_i} \to
\Ann_{C_i} x \to H^{i}_\mm(R)/x H^{i}_\mm(R) \to E^{t_i}/xE^{t_i} = 0.
\label{LCeqn1}
\end{align}
We also have, from the sequence $0 \to R \overset x \to R \to R/xR \to 0$,
the long exact sequence in local cohomology, which gives
\begin{align}
0 \to H^{i}_\mm(R)/x H^{i}_\mm(R) \to H^i_\mm(R/xR) \to
\Ann_{H^{i+1}_\mm(R)} x \to 0.  \label{LCeqn2}
\end{align}
From (\ref{LCeqn1}) and (\ref{LCeqn2}) we see that it suffices to bound
$\lambda\left(H^{i}_\mm(R)/x H^{i}_\mm(R)\right)$ and
$\lambda\left(\Ann_{H^{i+1}_\mm(R)} x\right)$ for $i \le d-2$,
and hence to bound
$\lambda(\Ann_E x)$ depending only on $\r$ (using the fact
that each $C_i$ embeds in a finite direct sum of $E$'s).
By duality, $\lambda\left(\Ann_E x\right)
= \lambda\left(R/(P+xR)\right)$.

Thus the problem reduces to showing that if $(S, \mm)$ is a one dimensional
complete domain and $x \in \mm - \mm^\r$ then there is a bound on
$\lambda(S/xS)$ which depends only on $\r$.  Let $T$ be the integral closure
of $S$, and set $N$ to be the degree of the extension of fraction fields.
Then $\lambda_S(S/xS) \le N \lambda_T(T/xT)$.  By Rees's strong valuation
theorem, there is an integer $k$ such that $\ord_T(x) \le \ord_S(x) + k$
\cite{Rees}.  This shows that $\lambda_S(S/xS) \le N(\r+k)$.
\end{proof}

\begin{thm} \label{CM1}
Let $(R,\mm,k)$ be a formally unmixed local ring of dimension $d$
such that $\A(\hat{R})$ is a prime ideal of dimension one in the
completion $\hat{R}$ of $R$. Then $R$ has a uniform bound on
relation type of parameter ideals.
\end{thm}

\begin{proof} By \cite[Lemma 4.1]{Wang1}, since our hypotheses pass to the completion, we may assume that $R$ is complete. By \cite[Lemma 2.2]{Wang:1997}, we may also assume that $H^0_\mm(R) = 0$. Observe further that we can assume $k$ is infinite. Indeed, let $S=R[y]_{m R[y]}$. Then the residue field of $S$ is infinite. Since $R \hookrightarrow S$ is smooth, for $c \in R$,  $R_c$ is CM if and only if $S_c$
is CM. Thus, $\A(R)S \subseteq \A(S)$. Moreover, for any prime
ideal $P \subseteq R$, $PS$ is prime. This implies that if $\A(S)$ properly
contains $\A(R)S$ then it is primary to the maximal ideal of $S$,
and so $S$ has f.l.c. It then follows that $R$ has f.l.c., a contradiction.
Hence, $\A(R)S = \A(S)$. We may pass from $R$ to $S$ and assume that $k$ is infinite.

Let $I=(x_1,\dots,x_d)$ be a parameter ideal in $R$. We may pick $x_d, x_{d-1}, \dots, x_1$ to form a superficial sequence. If $d = 2$, the theorem is true by \cite{Wang1}.

Suppose $d \ge 3$. Since $\dim R/\A(R)=1$, in $R/\A(R)$, $\mm + \A(R)$ has a principal minimal reduction $\overline{w}$. Let $\r$ be the reduction number of $\mm + \A(R)$ in $R/\A(R)$, and let $w$ be a representative of $\overline{w}$ in $R$. Suppose $y \in R$. Then, there exists an integer $t_y \ge 0$ such that $\overline{y} \in \overline{\mm}^{t_y} - \overline{\mm}^{t_y-1}$ in $R/\A(R)$. Observe that if $t_y \ge \r+1$, then we can write $y = y'w^l+ \alpha$ where $y' \in \mm^\r, l = t_y - \r$ and $\alpha \in \A(R)$. If $t_y \le \r$, we can write $y$ in the same form $y = y'w^l + \alpha$ by letting $y' = y, l = 0$ and $\alpha = 0$.

By replacing $x_d$ by a generic combination of $x_1, \dots, x_{d-1}$, we may first assume that $t_{x_d} = \min \{ t_{x_i} ~\big|~ 1 \le i \le d \}$ (since being superficial is an open condition, after replacing $x_d$ by a generic combination of $x_1, \dots, x_{d-1}$, the sequence $x_d, x_{d-1}, \dots, x_1$ is still superficial). We will now use Lemma \ref{rt-equality} to modify our parameter ideal as follows. If $t_{x_d} \ge \r+1$ then we write $x_d=x_d'w^t+\alpha$ where $x_d' \in \mm^\r, t = t_{x_d} - \r$ and $\alpha \in \A(R)$. We may choose $w, x_d'$ and $\alpha$ such that $(x_1, \dots, x_{d-1}, w)$ and $(x_1, \dots, x_{d-1}, x_d')$ are s.o.p. 's. That is, $(x_1, \dots, x_{d-1}, x_d'w^t)$ is a s.o.p. Let $y_i = x_i$ for $1 \le i \le d-1$, and $y_d = x_d'w^t$. It follows from Lemma \ref{rt-equality} that
$$ \rt (x_1, \dots, x_d) = \rt (y_1, \dots, y_d).$$
If $t_{x_d} \ge \r$, then we let $(y_1, \dots, y_d) = (x_1, \dots, x_d)$ and $x_d' = x_d$. The theorem will be proved if we can show that $\rt (y_1,y_2,\dots, y_d)$ is uniformly bounded. By Lemma \ref{superficial}, we may assume that $y_d, \dots, y_1$ form a superficial sequence.

It follows from Lemma \ref{CM1cohombound} that for $2 \le q \le d$ there exists a uniform bound $B_q$ (depending only on $q$ and $\r$) such that for $z_{q-1}, \dots, z_{d-1} \in R$ so that $(z_q, \dots, z_{d-1}, x_d')$ is part of a s.o.p., we have
$\lambda \big( H^0_\mm(R/(z_{q-1}, \dots, z_{d-1}, x_d')) \big) < B_q$. For $q = 2, \dots, d$, let
$$L_q = \max \{ B_q, B_{q+1} \}.$$
By considering the exact sequence
$0 \rightarrow R \stackrel{\gamma}{\rightarrow} R \rightarrow R/(\gamma) \rightarrow 0,$
for any non-zero-divisor $\gamma \in \A(R)$ which is part of a s.o.p., it also follows from Lemma \ref{CM1cohombound} that there exists a uniform bound $L_{d+1}$ not depending on $\gamma$ such that
$$\lambda\big( H^0_\mm(R/(\gamma)) \big) < L_{d+1}.$$
Let $L = (L_2-1) + \dots + (L_{d+1}-1)+1$. Let $M(d,k,l)$ be the Ramsey number determined by Lemma \ref{ramseynumber}. Similar to what was done in Theorem \ref{dparameters}, we recursively construct the following sequence of finite numbers: let $K_1 = M(d,1,L) = N_1$; and for $i \ge 2$, let $M_i = \sum_{l=0}^{K_{i-1}} {l+i-2 \choose i-2}$, $N_i = M(d, 2K_{i-1}, M_i(L-1)+1)$ and $K_i = 2K_{i-1}+N_i$.

Our proof now proceeds along a very similar line of argument (with some modification at the end) as in Theorem \ref{dparameters}. Use graded reverse lex monomial ordering with $y_1 > y_2 > \dots > y_d$ and $T_1 > T_2 > \dots > T_d$, and consider an arbitrary relation in $y_1, \dots, y_d$ of degree $N$
\begin{align}
F(T_1, \dots, T_d) = \sum_{|\n| = N} r_\n \T^\n = 0. \label{CM1start}
\end{align}
As in Theorem \ref{dparameters}, it suffices to show that for any $1 \le j \le d-1$, if $F$ contains a term $r_\m \T^\m$ (say $\m = (m_1, \dots, m_d)$) with $\ord_{< j+1}(\T^\m) > K_j$ then we can write $F(T_1, \dots, T_d) = H(T_1, \dots, T_d)G(T_1, \dots, T_d) + F'(T_1, \dots, T_d)$, where $H$ is a monomial divisible by $\prod_{i > j} T_i^{m_i}$, $G$ and $F'$ both provide relations in $x_1, \dots, x_d$, the leading term $r_\k \T^\k$ of $G$ satisfies the condition $r_\k \T^\k \big| r_\m \T^\m, \ord_{< j+1} (\T^\k) \le K_j$ and $\ord_{j+1}(\T^\k) = \dots = \ord_d(\T^\k) = 0$ (in particular, the degree of $G$ is bounded by $K_j$), and all terms $r_{\n_1}\T^{\n_1}$ of $F'$ with $\ord_{< j+1}(\T^{\n_1}) > K_j$ are smaller than $r_\m \T^\m$.

By contradiction, suppose our assertion is not true. As before, let $j$ be the smallest index for which there is a relation $F$ contradicting our assertion. Suppose $r_{\n_0} \T^{\n_0}$ is the largest term of $F$ for which $\ord_{< j+1}(\T^{\n_0}) > K_j$. We shall pick $F$ such that $r_{\n_0} \T^{\n_0}$ is smallest possible. Let $\n_0 = (n_{01}, \dots, n_{0d})$. For simplicity, we write $r$ for $r_{\n_0}$. We shall derive a contradiction.

From the choice of $j$, we may assume that $K = \ord_{< j}(\T^{\n_0}) \le K_{j-1}$. Again, we first observe that
\begin{align}
n_{0j} > K_j - K_{j-1} = K_{j-1}+N_j. \label{CM1n0}
\end{align}
Let $B = (y_{j+1}, \dots, y_d)$. We proceed along the same line of argument used in going from (\ref{gCMn0}) to (\ref{red1}), with the exception that instead of having $\A(R)$ we now have $\A(R)+(y_d)$ being $\mm$-primary. Hence, $y_j^{q_j'} - cy_d \in \A(R)$ for some positive integer $q_j'$ and $c \in R$, and so
$$(y_j^{q_j'} - cy_d) \Big(r \sfrac{\y^{\n_0}}{\prod_{i > j}y_i^{n_{0i}}} + \sum_{\n \in \J, \n < \n_0, ~\prod_{i > j}y_i^{n_{0i}} \big| \y^\n} r_{\n} \sfrac{\y^\n}{\prod_{i > j}y_i^{n_{0i}}}\Big) \in B.$$
Since $y_d \in B$, we get an equality similar to (\ref{red1})
$$y_j^{q_j'} \Big(r \sfrac{\y^{\n_0}}{\prod_{i > j}y_i^{n_{0i}}} + \sum_{\n \in \J, \n < \n_0, ~\prod_{i > j}y_i^{n_{0i}} \big| \y^\n} r_{\n} \sfrac{\y^\n}{\prod_{i > j}y_i^{n_{0i}}}\Big) \in B.$$
By a similar argument as in Theorem \ref{dparameters} again (as to get (\ref{neweq0}) and the condition in the next paragraph), we obtain a new relation in $y_1, \dots, y_d$,
\begin{align*}
Q(T_1, \dots, T_d) & = r T_1^{n_{01}} \dots T_{j-1}^{n_{0(j-1)}} T_j^{n_{0j}+q_j} + \sum_{|\m| = q_j+K+n_{0j}, \m < (n_{01}, \dots, n_{0j}+q_j, 0, \dots, 0)} u_\m \T^\m \\
& = r \T^{\m_0} + \sum_{\n < \m_0} u_\m \T^\m,
\end{align*}
where $\m_0 = (m_{01}, \dots, m_{0d}) = (n_{01}, \dots, n_{0(j-1)}, n_{0j}+q_j, 0, \dots, 0)$, and each term $u_\m \T^\m$ in $Q$ satisfies $\ord_{< j}(\T^\m) \le K_{j-1}$.

Let $\gamma \in \A(R)$ be a homology multiplier in $R$ such that $(\gamma, y_1, \dots, \widehat{y_j}, \dots, y_d)$ is a s.o.p. We also pick $\gamma$ such that $y_d$ is not in any associated primes of $\gamma$. Since $\gamma$ is part of a s.o.p. and $H^0_\mm(R) = 0$, $\gamma$ is a non-zero-divisor. Recall that $K = \ord_{< j}(\T^{\n_0}) = \ord_{< j}(\T^{\m_0}) \le K_{j-1}$. For each $1 \le p \le m_{0j}-K_{j-1}+1$, let $F_p(T_1, \dots, T_d)$ be the sum of all terms in $Q$ that are divisible by $T_j^{m_{0j}-K_{j-1}-p+1}$, and let
$$G_p(T_1, \dots, T_d) = \sfrac{F_p(T_1, \dots, T_d)}{T_j^{m_{0j}-K_{j-1}-p+1}} = r T_1^{m_{01}} \dots T_{j-1}^{m_{0(j-1)}} T_j^{K_{j-1}+p-1} + \ \text{smaller terms}.$$
We can continue in the same line of argument as in Theorem \ref{dparameters} (up to (\ref{gCMsteps})) to get a system of equalities which is similar to (\ref{gCMsteps})
\begin{eqnarray}
\left\{ \begin{array}{rccl}
\gamma G_1 & = & s_1 \y^{\m_1} + \sum_{\m < \m_1} s_\m \y^\m, \\
\gamma G_2 & = & s_2 \y^{\m_2} + \sum_{\m < \m_2} s_\m \y^\m, \\
& \vdots & & \\
\gamma G_p & = & s_p \y^{\m_p} + \sum_{\m < \m_p} s_\m \y^\m, \\
& \vdots & & \\
\gamma G_{m_{0j}-K_{j-1}+1} & = & s_{m_{0j}-K_{j-1}+1} \y^{\m_{m_{0j}-K_{j-1}+1}} + \sum_{\m < \m_{m_{0j}-K_{j-1}+1}} s_\m \y^\m,
\end{array} \right. \label{CM1steps}
\end{eqnarray}
where $\y^\m_p \in C_p = (y_1, \dots, \widehat{y_j}, \dots, y_d)^{K_{j-1}}B^{K+p}, |\m_p| = K_{j-1}+K+p$, and $G_p$ denotes $G_p(y_1, \dots, y_d)$, for $1 \le p \le m_{0j}-K_{j-1}+1$. Let $H_p = H_p(y_1, \dots, y_d)$ be the right hand side of the $p$-th equality in (\ref{CM1steps}).
As before, among all possible system of the form (\ref{CM1steps}) associated to the relation $Q(T_1, \dots, T_d)$, we choose one such that all the leading terms on the right hand side are minimal.

The following claim is similar to Claim \ref{gCMclaim}.

\begin{claim} \label{CM1claim}
There exists an integer $p \le N_j$ such that $H_p = 0$.
\end{claim}

\noindent{\sc Proof of Claim.} By (\ref{CM1n0}), we have $m_{0j}-K_{j-1}+1 \ge N_j+1$. By contradiction, suppose the assertion is false. That is, $s_1, \dots, s_{N_j}$ are all non-zero. Fix an integer $1 \le p \le m_{0j}-K_{j-1}+1$, and suppose the $y_i$'s that appear in $\y^{\m_p}$ are in $\{ y_i ~|~ i \ge h\} \backslash \{y_j\}$ (and $h$ is chosen to be the largest integer with this property).

If $h < d$, we have
\begin{align*}
s_p & \in \Big( \gamma, \{ \y^\m \in C_p ~ \big|~ \m < \m_p \} \Big) : \y^{\m_p} \\
& = \Big( \gamma, \{ \x^\m \in C_p ~|~ \ord_i(\y^\m) > \ord_i(\y^{\m_p}) \ \text{for some} \ i \ge h+1 \} \Big) : \y^{\m_p}.
\end{align*}
Choose $\alpha \in \A(R)$ such that the images of $(\gamma, \{ y_{h+1}, \dots, y_{d-1} \} \backslash \{ y_j \}, \alpha)$ form a part of a s.o.p.  in $R/(x_d')$ (this is possible because $h \not= j$ and $d \ge 3$). Then,
$$ \alpha s_p \in \big(\gamma, \{ y_{h+1}, \dots, y_d\} \backslash \{y_j\}\big) \subseteq \big( \gamma, \{ y_{h+1}, \dots, y_{d-1} \} \backslash \{ y_j \}, x_d' \big), $$
i.e.
$$ s_p \in \big(\gamma, \{ y_{h+1}, \dots, y_{d-1} \} \backslash \{y_j\}, x_d' \big) : \alpha.$$
Let $\tilde{\mm}, \tilde{s_p}, \tilde{\gamma}, \tilde{y_i} (i = 1, \dots, d),$ and $\tilde{\alpha}$ be the images of $\mm, s_p, \gamma, y_i (i=1, \dots, d)$, and $\alpha$ in $R/(x_d')$. Since $x_d'$ is part of a s.o.p., $R/(x_d')$ is a generalized Cohen-Macaulay ring. Since $(\tilde{\gamma}, \{ \widetilde{y_{h+1}}, \dots, \widetilde{y_{d-1}} \} \backslash \{ \tilde{y_j} \}, \tilde{\alpha})$ is part of a s.o.p., $\tilde{\alpha}$ avoids all associated primes of $(\tilde{\gamma}, \{ \widetilde{y_{h+1}}, \dots, \widetilde{y_{d-1}} \} \backslash \{ \tilde{y_j} \})$ except possibly $\tilde{\mm}$. Thus, we have
$$ \tilde{s_p} \in H^0_{\tilde{\mm}}\big( R\big/ \big(\tilde{\gamma}, \{ \widetilde{y_{h+1}}, \dots, \widetilde{y_{d-1}} \} \backslash \{ \tilde{y_j} \} \big) \big), $$
whence
$$ \overline{s_p} \in H^0_\mm\big(R\big/\big(\gamma, \{ y_{h+1}, \dots, y_{d-1} \} \backslash \{y_j\}, x_d' \big) \big), $$
here $\overline{s_p}$ denotes the image of $s_p$ in $R\big/\big(\gamma, \{ y_{h+1}, \dots, y_{d-1} \} \backslash \{y_j\}, x_d' \big)$.

If $h = d$, then the $p$-th equality of (\ref{CM1steps}) is
\begin{align}
\gamma G_p = s_p y_d^{m_p}, \label{CM1extra}
\end{align}
where $\m_p = (0, \dots, 0, m_p)$. Since $\gamma$ was chosen such that $y_d$ is not in any associated primes of $\gamma$, we must have
$$\overline{s_p} \in H^0_\mm(R/(\gamma)).$$

Consider the sequence $\{ \m_1, \m_2, \dots, \m_{N_j} \}$. Since $N_j = M(d,2K_{j-1},M_j(L-1)+1) \ge M(d,K_{j-1}+K,M_j(L-1)+1)$, it follows from Lemma \ref{ramseynumber} that there exists an increasing subsequence of length $V = M_j(L-1)+1$, $\m_{j_1} \preceq \m_{j_2} \preceq \dots \preceq \m_{j_V}$, with $1 \le j_1 < j_2 < \dots < j_V \le N_j \le m_{0j}-K_{j-1}$. Since $C_p = (y_1, \dots, \widehat{y_j}, \dots, y_d)^{K_{j-1}}B^{K+p}$, we have $\ord_{< j}(\y^{\m_p}) \le K_{j-1}$ for any $1 \le p \le m_{0j}-K_{j-1}+1$. Thus, since $M_j$ is the number of $(j-1)$-tuples of non-negative integers whose sum is at most $K_{j-1}$, we can choose from the sequence $\m_{j_1} \preceq \m_{j_2} \preceq \dots \preceq \m_{j_V}$ a subsequence of length $L$, $\m_{i_1} \preceq \m_{i_2} \preceq \dots \preceq \m_{i_L}$ such that they have the same $j-1$ first components, i.e. for any $1 \le a, b \le L$ we have $(\ord_1(\y^{\m_{i_a}}), \dots, \ord_{j-1}(\y^{\m_{i_a}})) = (\ord_1(\y^{\m_{i_b}}), \dots, \ord_{j-1}(\y^{\m_{i_b}}))$. Since $L = (L_2-1) + \dots + (L_{d+1}-1) + 1$, we can choose from this sequence a subsequence $\m_{p_1} \preceq \m_{p_2} \preceq \dots \preceq \m_{p_{L_q}}$, for some $2 \le q \le d+1$, such that if $q \le d$ then
$$\overline{s_{p_t}} \in H^0_\mm\big(R\big/\big(\gamma, \{ y_{q+1}, \dots, y_{d-1} \} \backslash \{y_j\}, x_d' \big)\big),$$
for all $t = 1, \dots, L_q$, and if $q = d+1$ then
$$\overline{s_{p_t}} \in H^0_\mm(R/(\gamma)),$$
for all $t = 1, \dots, L_q.$ Notice that when $q = d+1$, as in (\ref{CM1extra}), for all $t = 1, \dots, L_q$, the $p_t$-th equality of (\ref{CM1steps}) is
$$\gamma G_{p_t} = s_{p_t} y_d^{m_{p_t}},$$
where $\m_{p_t} = (0, \dots, 0, m_{p_t})$.

If $q = d+1$, it follows from the choice of $L_q$ that we can write
$$s_{p_{L_q}} = \sum_{i=1}^{L_q-1}b_is_{p_i} + a \gamma.$$
Substituting this into the equality $\gamma G_{p_{L_q}} = H_{p_{L_q}}$, we get
\begin{align}
\gamma \Big( G_{p_{L_q}} - ay_d^{m_{p_{L_q}}} \Big) = \sum_{i=1}^{L_q-1} \Big(b_i y_d^{m_{p_{L_q}} - m_{p_i}}\Big) s_{p_i} y_d^{m_{p_i}}. \label{CM1extrasim}
\end{align}
We can now use the $p_1$-th, $\dots$, $p_{L_q-1}$-th equalities in (\ref{CM1steps}) to simplify the right hand side of (\ref{CM1extrasim}) as we did in the $2$-dimensional case (Theorem \ref{2parameters}), bringing all terms with $\gamma$ to the left hand side and absorbing extra powers into the coefficient if necessary, to get
$$\gamma \Big(G_{p_{L_q}} + \ \text{terms smaller than} \ r y_1^{m_{01}} \dots y_{j-1}^{m_{0(j-1)}} y_j^{K_{j-1}+p_{L_q}-1}\Big) = 0.$$
This contradicts the fact that the right hand side of (\ref{CM1steps}) was chosen to be minimal.

Consider the case $q \le d$. By the choice of $L_q$, we can write
$$ s_{p_{L_q}} = \sum_{q < i < d, i \not= j} a_i y_i + a_j \gamma + a_d x_d' + \sum_{i=1}^{L_q-1}b_is_{p_i}. $$
Substitute this into the equality $\gamma G_{p_{L_q}} = H_{p_{L_q}}$, we get
\begin{align}
\gamma (G_{p_{L_q}} - a_j\y^{\m_{p_{L_q}}}) & = \sum_{q < i < d, i \not= j} a_iy_i \y^{\m_{p_{L_q}}} + a_d x_d' \y^{\m_{p_{L_q}}} + \sum_{i=1}^{L_q-1}b_is_{p_i}\y^{\m_{p_{L_q}}} + \sum_{\m < \m_{p_{L_q}}} s_\m \y^\m \nonumber \\
& = \sum_{q < i < d, i \not= j} (a_iy_q) \sfrac{\y^{\m_{p_{L_q}}}y_i}{y_q} + \sum_{i=1}^{L_q-1} \Big(b_i\sfrac{\y^{\m_{p_{L_q}}}}{\y^{\m_{p_i}}}\Big)s_{p_i}\y^{\m_{p_i}} + a_d x_d' \y^{\m_{p_{L_q}}} \nonumber \\
& \hspace{15ex} + \sum_{\m < \m_{p_{L_q}}} s_\m \y^\m. \label{CM1simplifying}
\end{align}
Since $\y^{\m_{p_{L_q}}} \in C_{p_{L_q}}$, (after absorbing extra powers into the coefficient if necessary) $a_j\y^{\m_{p_{L_q}}} < r y_1^{m_{01}} \dots y_{j-1}^{m_{0(j-1)}} y_j^{K_{j-1}+p_{L_q}-1}$. Also, for $q < i < d$,  $\sfrac{\y^{\m_{p_{L_q}}}y_i}{y_q} < \y^{\m_{p_{L_q}}}$. Furthermore, since $\y^{\m_{p_{L_q}}}$ and $\y^{\m_{p_i}}$ have the same $j-1$ first components, we have $\sfrac{\y^{\m_{p_{L_q}}}}{\y^{\m_{p_i}}} \in (y_{j+1}, \dots, y_d)$, which then implies that all terms of  $(b_i\sfrac{\y^{\m_{p_{L_q}}}}{\y^{\m_{p_i}}}\Big)G_{p_i}$ are smaller than $r y_1^{m_{01}} \dots y_{j-1}^{m_{0(j-1)}} y_j^{K_{j-1}+p_{L_q}-1}$. Thus, after using (\ref{CM1steps}) to replace $s_{p_i}\y^{\m_{p_i}}$ by $\big[ \gamma G_{p_i} - \sum_{|\m| = |\m_{p_i}|, \m < \m_{p_i}} s_\m \y^\m \big]$ in (\ref{CM1simplifying}), we get a relation
\begin{eqnarray}
\gamma\Big(G_{p_{L_q}} + \ \text{terms smaller than} \ r y_1^{m_{01}} \dots y_{j-1}^{m_{0(j-1)}} y_j^{K_{j-1}+p_{L_q}-1}\Big) = H_{p_{L_q}}' + a_d x_d' \y^{\m_{p_{L_q}}}, \label{CM1newrel}
\end{eqnarray}
where the leading term of $H_{p_{L_q}}'$ is strictly smaller than $s_{p_{L_q}}\x^{\m_{p_{L_q}}}$, which is the leading term of $H_{p_{L_q}}$.

If $x_d' = x_d$ (i.e. $x_d$ is in a small power of $\mm$), then $a_dx_d' \y^{\m_{p_{L_q}}} = (a_d y_q) \sfrac{\y^{\m_{p_{L_q}}}x_d'}{y_q}$. Thus, by absorbing $y_q$ into the coefficient, $a_dx_d' \y^{\m_{p_{L_q}}}$ gives a term strictly smaller than $\y^{\m_{p_{L_q}}}$. We can rename $H_{p_{L_q}}' + a_d x_d' \y^{\m_{p_{L_q}}}$ in (\ref{CM1newrel}) as $H_{p_{L_q}}'$. This contradicts the way (\ref{CM1steps}) was chosen.

Suppose now that $y_d = x_d'w^t$. Let us write $y_q = x_q'w^t + \beta$ for some $\beta \in \A(R)$ (this is possible since $t_{x_q} \ge t_{x_d}$). Write $\y^{\m_{p_{L_q}}} = \prod_{i \ge q, i \not= j} y_i^{l_i}$. Then, since all terms of $H_{p_{L_q}}'$ are smaller than $\y^{\m_{p_{L_q}}}$, it follows from (\ref{CM1newrel}) that
\begin{align*}
a_d & \in \big( \gamma, \{ y_i^{l_i+1} ~|~ i > q, i \not= j\} \big) : x_d' \y^{\m_{p_{L_q}}} \\
& \subseteq \big( \gamma, \{ y_i^{l_i+1} ~|~ q < i < d, i \not= j\}, (w^t)^{l_d+1} \big) : x_d' \y^{\m_{p_{L_q}}}.
\end{align*}
Therefore,
\begin{align}
\beta a_d \in \big(\gamma, \{ y_i ~|~ q < i < d, i \not= j \}, w^t \big).\label{CM1newHM}
\end{align}
We also have
\begin{align}
a_dx_d'\y^{\m_{p_{L_q}}} & = a_dx_d'(x_q'w^t+\beta)y_q^{l_q-1} \big( \prod_{q < i < d, i \not= j} y_i^{l_i} \big) ~ (x_d'w^t)^{l_d} \nonumber \\
& = a_dx_d'x_q'w^t y_q^{l_q-1} \big( \prod_{q < i < d, i \not= j} y_i^{l_i} \big) ~ (x_d'w^t)^{l_d} + \beta a_d x_d' y_q^{l_q-1} \big( \prod_{q < i < d, i \not= j} y_i^{l_i} \big) ~ (x_d'w^t)^{l_d} \nonumber \\
& = (a_dx_q')(x_d'w^t)^{l_d+1} y_q^{l_q-1} \big( \prod_{q < i < d, i \not= j} y_i^{l_i} \big) + \beta a_d x_d' y_q^{l_q-1} \big( \prod_{q < i < d, i \not= j} y_i^{l_i} \big) ~ (x_d'w^t)^{l_d} \nonumber \\
& = (a_dx_q') y_q^{l_q-1} \big( \prod_{q < i < d, i \not= j} y_i^{l_i} \big) ~ y_d^{l_d+1} + \beta a_d x_d' y_q^{l_q-1} \big( \prod_{q < i < d, i \not= j} y_i^{l_i} \big) ~ (x_d'w^t)^{l_d}. \label{CM1moresim}
\end{align}
By (\ref{CM1newHM}), we can write
$$\beta a_d = c_j \gamma + \sum_{q < u < d, u \not= j} c_u y_u + c_d w^t.$$
Substituting this into (\ref{CM1moresim}), we get
\begin{align}
a_dx_d'\y^{\m_{p_{L_q}}} & = (a_dx_q') y_q^{l_q-1} \big( \prod_{q < i < d, i \not= j} y_i^{l_i} \big) ~ y_d^{l_d+1} + c_j \gamma x_d' y_q^{l_q-1} \big( \prod_{q < i < d, i \not= j} y_i^{l_i} \big) ~ (x_d'w^t)^{l_d} \nonumber \\
& \quad + \sum_{q < u < d, i \not= j} (c_ux_d') y_u y_q^{l_q-1} \big( \prod_{q < i < d, i \not= j} y_i^{l_i} \big) ~ (x_d'w^t)^{l_d} \nonumber \\
& \quad + c_d (x_d'w^t) y_q^{l_q-1} \big( \prod_{q < i < d, i \not= j} y_i^{l_i} \big) ~ (x_d'w^t)^{l_d} \nonumber \\
& = (a_dx_q') y_q^{l_q-1} \big( \prod_{q < i < d, i \not= j} y_i^{l_i} \big) ~ y_d^{l_d+1} + c_j \gamma x_d' y_q^{l_q-1} \big( \prod_{q < i < d, i \not= j} y_i^{l_i} \big) ~ y_d^{l_d} \nonumber \\
& \quad + \sum_{q < u < d, i \not= j} (c_ux_d') y_u y_q^{l_q-1} \big( \prod_{q < i < d, i \not= j} y_i^{l_i} \big) ~ y_d^{l_d} + c_d y_q^{l_q-1} \big( \prod_{q < i < d, i \not= j} y_i^{l_i} \big) ~ y_d^{l_d+1}. \label{CM1lastterm}
\end{align}
Observe that since $\y^{\m_{p_{L_q}}} \in (y_1, \dots, \widehat{y_j}, \dots, y_d)^{K_{j-1}}B^{K+p_{L_q}}$, $y_q^{l_q-1} \big( \prod_{q < i < d, i \not= j} y_i^{l_i} \big) y_d^{l_d} = \sfrac{\y^{\m_{p_{L_q}}}}{y_q} \in (y_1, \dots, \widehat{y_j}, \dots, y_d)^{K_{j-1}}B^{K+p_{L_q}-1},$ i.e.
$$(c_jx_d') y_q^{l_q-1} \big( \prod_{q < i < d, i \not= j} y_i^{l_i} \big) y_d^{l_d} < r y_1^{m_{01}} \dots y_{j-1}^{m_{0(j-1)}} y_j^{K_{j-1}+p_{L_q}-1}.$$
Hence, substituting (\ref{CM1lastterm}) into (\ref{CM1newrel}), bringing terms with $\gamma$ to the left hand side, we get
\begin{align}
\gamma \Big(G_{p_{L_q}} + \ \text{terms smaller than} \ r y_1^{m_{01}} \dots y_{j-1}^{m_{0(j-1)}} y_j^{K_{j-1}+p_{L_q}-1}\Big) = H_{p_{L_q}}'', \label{CM1lastrel}
\end{align}
where the leading terms of $H_{p_{L_q}}''$ is strictly smaller than $\y^{\m_{p_{L_q}}}$ which is the leading term of $H_{p_{L_q}}$.
This contradicts the way (\ref{CM1steps}) was chosen. Our claim is proved. \qed

We have just shown that there must exist an integer $p \le N_j$ such that
$$ \gamma G_p = 0.$$
Thus, since $\gamma$ is a non-zero-divisor, we have
$G_p = 0.$
That is,
$$ r y_1^{m_{01}} \dots y_{j-1}^{m_{0(j-1)}} y_j^{K_{j-1}+p-1} + \ \text{smaller terms} \ = 0. $$
This gives a new relation in $y_1, \dots, y_d$, which by abusing language we also denote by $G_p$,
\begin{align*}
G_p(T_1, \dots, T_d) & = r T_1^{m_{01}} \dots T_{j-1}^{m_{0(j-1)}} T_j^{K_{j-1}+p-1} + \ \text{smaller terms} \\
& = r T_1^{n_{01}} \dots T_{j-1}^{n_{0(j-1)}} T_j^{K_{j-1}+p-1} + \ \text{smaller terms}.
\end{align*}
Once again, observe that $K_{j-1}+p-1 \le K_{j-1}+N_j-1 < n_{0j}$ by (\ref{CM1n0}). Therefore, we can write
$$F(T_1, \dots, T_d) = \big[T_j^{n_{0j}-K_{j-1}-p+1}\prod_{i > j} T_i^{n_{0i}}\big] G_p(T_1, \dots, T_d) + F'(T_1, \dots, T_d),$$
where $G_p(T_1, \dots, T_d)$ and $F'(T_1, \dots, T_d)$ are relations in $y_1, \dots, y_d$. Again, it is clear that every term $r_\n \T^\n$ of $F'$ with $\ord_{< j+1}(\T^\n) > K_j$ is smaller than $r_{\n_0} \T^{\n_0}$. We obtain a contradiction. Hence, the theorem is proved.
\end{proof}

Let $R$ be a ring of positive prime characteristic $p$. We denote
the $e$th power of the Frobenius endomorphism $f:R \to R$ sending
$x \mapsto x^p$ by $f^e$. For $q = p^e$, a power of $p$, and $I
\subseteq R$ we let $I^{[q]} = (i^q | i \in I)$.  The ideal $I
\subseteq R$ is called {\sl Frobenius closed} if whenever $x^q \in
I^{[q]}$ then $x \in I$.  The ring $R$ is called   {\sl $F$-pure}
if $f$ is a pure morphism and {\sl cyclically  $F$-pure} if  all
ideals of $R$ are Frobenius closed.  When $R$ is excellent these
conditions are equivalent \cite{Hoch77}.

Also, when $R = S/J$ is the image of a regular local ring
$(S,\nn)$ then Fedder has given a criterion for $F$-purity in
terms of $J$ \cite{Fedder}.  $R$ is $F$-pure if and only if
$J^{[p]}:_S J \not\subseteq \nn^{[p]}$.

When $(R,\mm)$ is an excellent local ring then $R$ is $F$-pure if
and only if $\widehat R$ is $F$-pure.  Moreover, in an $F$-pure
ring, the ideal $\A(R)$ is radical.  We may thus apply
Theorem~\ref{CM1} to obtain

\begin{cor} \label{Fpure}
Let $(R,\mm)$ be an complete local equidimensional
$F$-pure ring such
that the defining  ideal of the non-CM locus is a dimension one
prime ideal.  Then $R$ has a uniform bound on relation type of
parameter ideals.
\end{cor}


\end{document}